\documentclass[12pt]{amsart}
\usepackage{enumitem}
\usepackage{psfrag}
\usepackage{graphicx}
\usepackage{pinlabel}
\usepackage{graphicx}
\usepackage{amsmath}
\usepackage{amssymb}
\usepackage{amscd}
\usepackage{autobreak}
\usepackage{picinpar}
\usepackage{times}
\usepackage{pb-diagram}
\usepackage{graphicx}
\usepackage{wrapfig}
\usepackage{xspace}
\usepackage{hyperref}
\usepackage{url}
\usepackage{caption}
\usepackage{subcaption}
\usepackage{fancyhdr}
\usepackage{overpic}
\usepackage{color}
\usepackage[square,comma,sort&compress,numbers]{natbib}
\usepackage{latexsym}
\usepackage{mathrsfs}
\usepackage{amsfonts}
\usepackage{hyperref}
\usepackage{amsthm}
\usepackage{appendix}
\usepackage{pdfsync}
\usepackage{lineno}

\theoremstyle{definition}
\newtheorem{theorem}{Theorem}[section]
\newtheorem{lemma}[theorem]{Lemma}

\newtheorem{corollary}[theorem]{Corollary}
\newtheorem{definition}[theorem]{Definition}
\newtheorem{conjecture}[theorem]{Conjecture}
\newtheorem{remark}[theorem]{Remark}
\newtheorem{convention}[theorem]{Convention}
\newtheorem*{theorem*}{Theorem}

\graphicspath{{figures/}}   

\def\qed{\hfill{Q.E.D.}\smallskip}

\begin{document}

\title{\bf Maximal principles in discrete conformal geometry with application to the rigidity of infinite triangulations}
\author{Yanwen Luo, Xu Xu, Chao Zheng}

\date{\today}
\address{Department of Mathematics, Oklahoma State University, Stillwater, 74074, U.S.} \email{yanwen.luo@okstate.edu}

\address{School of Mathematics and Statistics, Wuhan University, Wuhan, 430072, P.R.China} \email{xuxu2@whu.edu.cn}

\address{School of Mathematics and Statistics, Wuhan University, Wuhan 430072, P.R. China} \email{czheng@whu.edu.cn}



\begin{abstract}
In this paper, maximum principles for Euclidean and hyperbolic discrete conformal structures on polyhedral surfaces are established.
These maximum principles unify and generalize the maximum principles for vertex scalings and different types of circle packings in the literature.
As an application of the hyperbolic discrete maximum principle, a discrete Schwarz-Ahlfors lemma is established.
As another application, an infinite rigidity theorem for small Delaunay triangulations of the hyperbolic plane is proved.
\end{abstract}

\maketitle

\tableofcontents

\section{Introduction}

The maximum principle serves as a ubiquitous tool in the study of various problems in analysis and geometry, including solutions to partial differential equations, minimal surfaces, harmonic maps and others.
In complex analysis, the maximum principle directly leads to the Schwarz lemma and Liouville's theorem,  two fundamental properties of holomorphic functions.
It can also be used to characterize the conformal automorphisms of the unit disk and the complex plane.

In the novel field of discrete conformal geometry, various forms of maximum principles have been established to prove fundamental results regarding the existence and uniqueness of discrete conformal maps \cite{Thurston, Gu1, Gu2, BL-E, BL-H}, as well as their convergence to smooth counterparts \cite{RS, HS1996, LSW, Bucking}.
In this paper, we explore the maximum principle for general discrete conformal structures of Euclidean and hyperbolic polyhedral surfaces, as proposed in \cite{Glickenstein JDG, G-T, ZGZLYG}.
This work generalizes previous maximum principles for Thurston's circle packings \cite{Thurston}, Luo's vertex scalings \cite{Luo CCM}, and Bowers-Stephenson's  inversive distance circle packings \cite{BS}.
As applications, we establish a discrete Schwarz-Ahlfors lemma and an infinite rigidity theorem for small Delaunay triangulations of the hyperbolic plane.

\subsection{A Euclidean discrete maximum principle}\label{Sec: 1.1}

Let $(S,\mathcal{T})$ be a triangulated surface, possibly with boundary, equipped with a triangulation $\mathcal{T}=\{V,E,F\}$,
where $V$, $E$, and $F$ represent the sets of vertices, edges, and faces, respectively.
Denote a vertex, an edge, and a face in the triangulation $\mathcal{T}$ by $v_i$, $v_iv_j$, and $\triangle v_iv_jv_k$, respectively.

A piecewise linear metric (abbreviated as PL metric) on $(S, \mathcal{T})$ is a function $l: E \rightarrow \mathbb{R}_{>0}$,
such that each face $\triangle v_iv_jv_k$ in $F$ is associated with a non-degenerate Euclidean triangle with edge lengths $l_{ij}$, $l_{ik}$, and $l_{jk}$.
For a PL metric $l: E \rightarrow \mathbb{R}_{>0}$ on $(S, \mathcal{T})$,
the combinatorial curvature is a map $K: V \rightarrow (-\infty, 2\pi)$,
which assigns to an interior vertex $v_i \in V$ the value $2\pi$ minus the sum of angles of Euclidean triangles at $v_i$,
and to a boundary vertex $v_i \in V$ the value $\pi$ minus the sum of angles at $v_i$.
An interior vertex $v$ is flat in a PL metric if $K(v)=0$.
A PL metric is flat if all interior vertices are flat.

\begin{definition}[\cite{G-T,XZ,ZGZLYG}]
\label{Def: IDCP}
Let $(S, \mathcal{T}, \varepsilon, \eta)$ be a weighted triangulated surface with the weights $\varepsilon: V \rightarrow \{-1,0,1\}$ and $\eta: E \rightarrow \mathbb{R}$ satisfying $\eta_{ij}=\eta_{ji}$.
A Euclidean discrete conformal structure on $(S, \mathcal{T}, \varepsilon, \eta)$ is defined by a map $f: V \rightarrow \mathbb{R}$, such that a PL metric $l: E \rightarrow \mathbb{R}_{>0}$ assigns the edge length $l_{ij}$ for each edge $v_iv_j \in E$ by
\begin{equation}\label{Eq: length-E}
l_{ij}=\sqrt{\varepsilon_ie^{2f_i}+\varepsilon_je^{2f_j}
+2\eta_{ij}e^{f_i+f_j}}.
\end{equation}
The function $f: V \to \mathbb{R}$ is referred to as a Euclidean label, and $\eta$ is the Euclidean weight on $(S, \mathcal{T}, \varepsilon, \eta)$.
Two PL metrics $(\varepsilon, \eta, l)$ and $(\bar{\varepsilon}, \bar{\eta}, \bar{l})$ on $(S, \mathcal{T})$ are said to be Euclidean discrete conformal equivalent if $\varepsilon=\bar{\varepsilon}$ and $\eta=\bar{\eta}$.
In this case, we define $w^e=\bar{f}-f$ and denote this relationship by $\bar{l}=w^e*l$, where $w^e$ is called a Euclidean discrete conformal factor.
\end{definition}

\begin{figure}[!ht]
\centering
\includegraphics[scale=1]{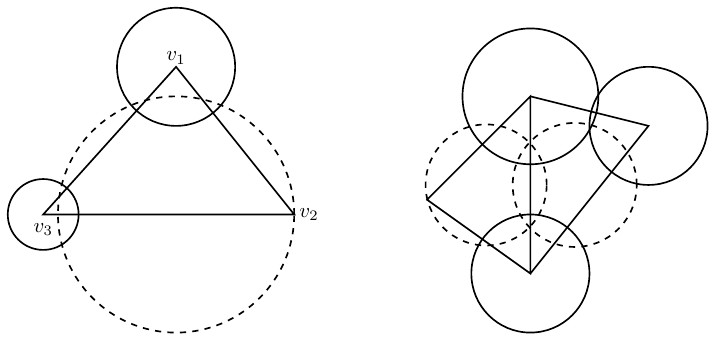}
\caption{A triangle with $\epsilon_2 = 0, \epsilon_1 = \epsilon_3=1$ and weighted Delaunay condition.}
\label{trianglequad}
\end{figure}

This definition unifies previous definitions of discrete conformal structures, including circle packings, vertex scaling, and inversive distance circle packings.
Specifically, let $r_i=e^{f_i}$ for all $i\in V$.
Then (\ref{Eq: length-E}) is equivalent to
\begin{equation}\label{Eq: length-E2}
l_{ij}=\sqrt{\varepsilon_ir_i^2+\varepsilon_jr_j^2+2\eta_{ij}r_ir_j}.
\end{equation}
If $\varepsilon\equiv0$, then (\ref{Eq: length-E2}) reduces to the vertex scaling of PL metric introduced by Luo \cite{Luo CCM}.
If $\varepsilon\equiv1$, then (\ref{Eq: length-E2}) reduces to the circle packings introduced in \cite{Thurston,BS}.

To see a unified geometric interpretation,
consider a non-degenerate Euclidean triangle $\triangle v_1v_2v_3$ induced by $(r_1,r_2,r_3)\in\mathbb{R}^3_{>0}$ via (\ref{Eq: length-E2}).
From a geometric perspective, if $\varepsilon\equiv1$, three circles with radii $r_1,r_2,r_3$ are located at the vertices of $\triangle v_1v_2v_3$, referred to as vertex-circles and denoted by $C_1,C_2,C_3$, respectively.
For any vertex $v_i$ with $\varepsilon_i=0$, the radius of the corresponding vertex-circle $C_i$ can be considered as $0$, implying that $C_i$ shrinks to a single point.  Note that this does not imply $r_i=0$, since $r_i=e^{f_i}>0$ by definition. Please refer to the triangle in Figure \ref{trianglequad}.

Next, we recall the definition of weighted Delaunay triangulations on $(S, \mathcal{T}, \varepsilon, \eta)$.
For a non-degenerate Euclidean triangle $\triangle v_1v_2v_3$  induced by $(r_1,r_2,r_3)\in\mathbb{R}^3_{>0}$ via (\ref{Eq: length-E2}), there exists a geometric center $c_{123}$ with the same power distances to the vertices $v_1,v_2,v_3$ (\cite{Glickenstein}, Proposition 7).
Here the power distance of a point $p$ to the vertex $v_i$ is defined to be $\pi_i(p)=|p-v_i|^2-\varepsilon_ir_i^2$.
The circle centered at $c_{123}$ with radius $\sqrt{\pi_i(c_{123})}$ is referred to as the face-circle of $\triangle v_1v_2v_3$, denoted by $C_{123}$.
It is straightforward to verify that the face-circle is orthogonal to the vertex-circle $C_i$ with $\varepsilon_i=1$ and passes through the vertex $v_j$ with $\varepsilon_j=0$.
Specifically, if $\varepsilon\equiv0$, then the face-circle $C_{123}$ is the circumscribed circle of $\triangle v_1v_2v_3$.

Let $h_{ij,k}$ be the signed distance from the geometric center $c_{123}$ to the edge $v_iv_j$, which is positive if $c_{123}$ lies on the same side of the line along $v_iv_j$ as $\triangle v_1v_2v_3$, and negative otherwise (or zero if $c_{ijk}$ is on the line along $v_iv_j$).
Please refer to Figure \ref{figure5} for an illustration of $h_{ij,k}$.
The common edge $v_1v_2$ of two adjacent non-degenerate Euclidean triangles $\triangle v_1v_2v_3$ and $\triangle v_1v_2v_4$ is called \textit{weighted Delaunay} in the PL metric $l$ if
\begin{equation}\label{Eq: weighted Delaunay}
h_{12,3}+h_{12,4}\geq 0.
\end{equation}
The weighted Delaunay condition (\ref{Eq: weighted Delaunay}) can be interpreted geometrically.
The edge $v_1v_2$ in the PL metric is weighted Delaunay if and only if, for $\varepsilon_4=0$, the vertex $v_4$ is not contained in the interior of the face-circle $C_{123}$, and for $\varepsilon_4=1$, the vertex-circle $C_4$ either does not intersect the face-circle $C_{123}$, or intersects it with an exterior angle of at most $\frac{\pi}{2}$.
Specifically, if $\varepsilon\equiv0$, the weighted Delaunay condition reduces to the classical Delaunay condition. See Figure \ref{trianglequad}.

The weighted triangulation $\mathcal{T}$ on $(S, \mathcal{T}, \varepsilon, \eta)$ with a PL metric is called weighted Delaunay if all interior edges are weighted Delaunay.
For simplicity, we refer to the PL metric $l$ or the labels $f$ ($r$ respectively) as weighted Delaunay.
It is important to note that (\ref{Eq: weighted Delaunay}) applies only to non-degenerate triangles.
In Subsection \ref{Sec: weighted Delaunay}, we will extend this definition to generalized triangles, i.e., some triangles may be degenerate.

\begin{figure}[!ht]
\centering
\includegraphics[scale=1.2]{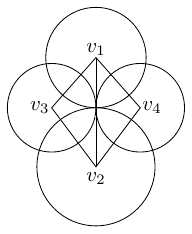}
\caption{The configuration of the four circles.}
\label{exception}
\end{figure}

In this paper, we focus on the weights $\varepsilon: V \rightarrow \{0,1\}$ and $\eta: E \rightarrow \mathbb{R}_{\geq 0}$ satisfying
\begin{equation}\label{Eq: SC1}
\varepsilon_s \varepsilon_t + \eta_{st} > 0, \quad \forall v_s v_t \in E.
\end{equation}
Unless otherwise stated, the weights $\varepsilon$ and $\eta$ throughout the paper will satisfy these conditions.
The weight pair $(\varepsilon, \eta)$ on $(S, \mathcal{T}, \varepsilon, \eta)$ is said to be \textit{regular} if there exists no pair of triangles $\triangle v_1v_2v_3$ and $\triangle v_1v_2v_4$ satisfying
\begin{equation*}
\varepsilon_1 = \varepsilon_2 = \varepsilon_3 = \varepsilon_4 = 1, \quad \eta_{12} = 1, \quad \eta_{13} = \eta_{23} = \eta_{14} = \eta_{24} = 0.
\end{equation*}
For the exceptional case, see Figure \ref{exception}.

Let $P_n$ be an $n$-sided star-shaped polygon with cyclically ordered boundary vertices $v_1,\dots,v_n$ (where $v_{n+1}=v_1$).
Let $v_0$ be an interior point of $P_n$, which induces a triangulation $\mathcal{T}$ of $P_n$ consisting of triangles $\triangle v_0v_iv_{i+1}$.
We have the following Euclidean maximal principle.

\begin{theorem}\label{Thm: MP-E}
Let $l,\bar{l}$ be two weighted Delaunay PL metrics on $(P_n,\mathcal{T},\varepsilon,\eta)$ with a regular weight pair $(\varepsilon,\eta)$, namely, $\bar{l}=w^e*l$.
If $K_0(l)\leq K_0(\bar{l})$,
then
\begin{gather*}
w^e_0\leq \max_{j \in\{ 1, 2, \cdots, n\}}w^e_j.
\end{gather*}
Furthermore, if $w^e_0 = \max_{j \in\{ 1, 2, \cdots, n\}}w^e_j$,
then $w^e$ is a constant on $P_n$.
\end{theorem}

\subsection{A hyperbolic discrete maximum principle and a discrete Schwarz-Ahlfors lemma}\label{Sec: 1.2}

A piecewise hyperbolic metric (abbreviated as PH metric) on $(S, \mathcal{T}, \varepsilon, \eta)$ is a function $l: E \rightarrow \mathbb{R}_{>0}$ such that each face in $F$ is associated with a non-degenerate hyperbolic triangle, analogous to the case of PL metrics.
The combinatorial curvature in the PH metric is defined in the same way as the Euclidean case, still denoted by $K$.

\begin{definition}[\cite{G-T,ZGZLYG}]
\label{Def: IDCP-H}
A hyperbolic discrete conformal structure on $(S, \mathcal{T}, \varepsilon, \eta)$ is defined by a map $f: V\rightarrow \mathbb{R}$ such that
a PH metric $l:E\rightarrow \mathbb{R}_{>0}$ assigns the edge length $l_{ij}$ for each edge $v_iv_j \in E$ by
\begin{equation}\label{Eq: length-H}
\cosh l_{ij}
=\sqrt{1+\varepsilon_ie^{2f_i}}
\sqrt{1+\varepsilon_je^{2f_j}}
+\eta_{ij}e^{f_i+f_j}.
\end{equation}
The function $f: V \to \mathbb{R}$ is referred to as the hyperbolic label, and $\eta$ is the hyperbolic weight on $(S, \mathcal{T}, \varepsilon, \eta)$.
\end{definition}

For any vertex $v_i\in V$, we set
\begin{equation}\label{Eq: u f}
u_i=\begin{cases}
f_i, & \text{if $\varepsilon_i=0$}, \\
\frac{1}{2} \ln \left|\frac{\sqrt{1+e^{2f_i}}-1}{\sqrt{1+e^{2f_i}}+1}\right|, & \text{if $\varepsilon_i=1$.}
\end{cases}
\end{equation}
Two PH metrics $(\varepsilon,\eta, l)$ and $(\bar\varepsilon,\bar\eta, \bar l)$ on $(S, \mathcal{T})$ are said to be hyperbolic discrete conformal equivalent if $\varepsilon=\bar\varepsilon$ and $\eta=\bar\eta$.
In this case, we define $w^h=\bar{u}-u$ and denote this relationship by $\bar l=w^h*l$, where $w^h$ is called a hyperbolic discrete conformal factor.

Similar to the Euclidean case, we have the following geometric explanation for the hyperbolic discrete conformal structures.
For a non-degenerate hyperbolic triangle $\triangle v_1v_2v_3$ induced by a hyperbolic discrete conformal structure in Definition \ref{Def: IDCP-H}, if $\varepsilon_i=1$, we set $\sinh r_i=e^{f_i}$ and thus $u_i=\ln \tanh \frac{r_i}{2}$.
Under this condition, there exists a hyperbolic vertex-circle centered at $v_i$ with radius $r_i$.
If $\varepsilon_i=0$, then we think the hyperbolic vertex-circle shrinks into the vertex $v_i$.

The definition of the weighted Delaunay condition in the hyperbolic case parallels that in the Euclidean case.
For a non-degenerate hyperbolic triangle $\triangle v_1v_2v_3$ induced by a hyperbolic discrete conformal structure in Definition \ref{Def: IDCP-H}, there is a geometric center $c_{123}$ similar to the Euclidean case. Please refer to Glickenstein-Thomas' work \cite{G-T} for more details. Note that the geometric center $c_{123}$ may be out of the hyperbolic space in this case.
By projecting the geometric center $c_{123}$ to the edges, one can define the signed distance $h_{ij,k}$ of the geometric center $c_{123}$
to the edge $v_iv_j$ similarly. The explicit expression of $h_{ij,k}$ in the hyperbolic case can be found in \cite{G-T,Xu MRL,Xu 1}.
If an edge $v_iv_j$ is shared by two non-degenerate hyperbolic triangles $\triangle v_iv_jv_k$ and $\triangle v_iv_jv_l$ generated by a hyperbolic discrete conformal structure in Definition \ref{Def: IDCP-H}, the edge $v_iv_j$ is said to be weighted Delaunay if
$$h_{ij,k}+h_{ij,l}\geq0.$$
The weighted Delaunay condition could also be represented by the conditions on face-circles in Figure \ref{trianglequad}.
The edge $v_1v_2$ in the PH metric is weighted Delaunay if and only if, for $\varepsilon_4=0$, the vertex $v_4$ is not contained in the interior of the face-circle $C_{123}$, and for $\varepsilon_4=1$, the vertex-circle $C_4$ either does not intersect the face-circle $C_{123}$, or intersects it with an exterior angle of at most $\frac{\pi}{2}$.
Note that in the Poincar\'{e} disk model of the hyperbolic plane $\mathbb{D}$, intersection angles of circles are the same as in the Euclidean background geometry, and
a hyperbolic circle coincides with a Euclidean circle, although the centers may be different.
Hence, an edge is weighted Delaunay in the PH metric if and only if it is weighted Delaunay in the induced PL metric.
Please refer to Section \ref{Sec: MP-H} for more details.

We have the following hyperbolic discrete maximum principle. 

\begin{theorem}\label{Thm: MP-H}
Let $l,\bar{l}$ be two weighted Delaunay PH metrics on $(P_n,\mathcal{T},\varepsilon,\eta)$ with a regular weight pair $(\varepsilon,\eta)$, namely, $\bar{l}=w^h*l$.
All circles or points with respect to $f,\bar f$ are contained in $\mathbb{D}$.
\begin{description}
  \item[(a)] If $K_0(l)\geq K_0(\bar{l})$ and $w^h_0=\bar{u}_0-u_0>0$, then
\begin{equation*}
w^h_0<\max_{j\in \{1, 2, \cdots, n\}}w^h_j.
\end{equation*}
  \item[(b)] If $K_0(l)\leq K_0(\bar{l})$ and $w_0^h=\bar{u}_0-u_0<0$, then
\begin{equation*}
w^h_0>\min_{j\in \{1, 2, \cdots, n\}}w^h_j.
\end{equation*}
\end{description}
\end{theorem}

By combining the formula $w^h = \bar{u} - u$ with the relationship between $u$ and $f$ in \eqref{Eq: u f}, Theorem \ref{Thm: MP-H} immediately implies the following discrete Schwarz-Ahlfors lemma. 
\begin{theorem}[Discrete Schwarz-Ahlfors lemma]\label{Thm: DSL2}
Let $(\varepsilon,\eta)$ be a regular pair of weights on $(M,\mathcal{T},\varepsilon,\eta)$, where $M\subseteq\mathbb{D}$ is a compact set with non-empty boundary.
Suppose $l$ and $\bar{l}$ are two weighted Delaunay PH metrics with hyperbolic labels $f$ and $\bar{f}$, respectively, satisfying all circles and points with respect to $f,\bar f$ are contained in $\mathbb{D}$.
\begin{description}
\item[(a)] If the combinatorial curvatures $K(l)\geq K(\bar{l})$ for all interior vertices, and $f \geq \bar{f}$ holds for every boundary vertex, then $f \geq \bar{f}$ holds for all vertices.
\item[(b)] If $K(l)\leq K(\bar{l})$ for all interior vertices, and $f \leq \bar{f}$ holds for every boundary vertex, then $f \leq \bar{f}$ holds for all vertices.
\end{description}
\end{theorem}
We will explain its relationship with the classical Schwarz-Ahlfors lemma after its proof in Section \ref{Sec: MP-H}.




\subsection{Rigidity of infinite triangulations of the plane} 

By applying a special case of the hyperbolic maximum principle in Theorem \ref{Thm: MP-H},
we establish a rigidity theorem for infinite geodesic triangulations on the hyperbolic disk $\mathbb{D}$.
If $\varepsilon \equiv 0$, then the edge length in (\ref{Eq: length-H}) reduces to $\cosh l_{ij}=1+\eta_{ij}e^{f_i+f_j} =1+\eta_{ij}e^{u_i+u_j}  $, which is equivalent to
\begin{equation*}
\sinh \frac{l_{ij}}{2}
=\sqrt{\frac{\eta_{ij}}{2}}e^{\frac{1}{2}(u_i+u_j)}.
\end{equation*}
Therefore, if two PH metrics $l$ and $\bar l$ are discrete conformal equivalent, then
\begin{equation}\label{Eq: DCE-H}
\sinh\frac{\bar l_{ij}}{2}
=e^{\frac{1}{2}(w^h_i+w^h_j)}\sinh \frac{l_{ij}}{2}.
\end{equation}
This corresponds to the vertex scaling of PH metrics introduced by Bobenko-Pinkall-Springborn \cite{BPS}.
The weighted Delaunay triangulations in this case coincide with the classical Delaunay triangulations.

\begin{theorem}\label{Thm: IR-H}
Let $l$ and $\bar{l}$ be two discrete conformal equivalent Delaunay PH metrics on $(\mathbb{D}, \mathcal{T})$, i.e., they are related by (\ref{Eq: DCE-H}).
If there exists a constant $\delta>0$ such that
\begin{enumerate}
\item[(a)] all inner angles of all triangles in $l$ and $\bar{l}$ are at least $\delta$,
\item[(b)] all edge lengths of all triangles in $l$ and $\bar{l}$ are less than $\delta^3/8192$,
\item[(c)] the combinatorial curvatures $K(l)\equiv K(\bar{l})\equiv0$ for all vertices,
\end{enumerate}
then $l=\bar{l}$, i.e., they are isometric in $\mathbb{D}$.
\end{theorem}


Theorem \ref{Thm: IR-H} can be regarded as a discrete analogue of the following well-known result in complex analysis: every conformal automorphism of the unit disk is a hyperbolic isometry.
Similar rigidity theorems have been established for circle packings \cite{RS, He1991, Sch, He, HS1993, HS1995, HS1996}, vertex scalings \cite{WGS, LSW, DGM, Dai-Wu}, and inversive distance circle packings \cite{LXZ}.
These rigidity results are essential for proving that the $K$-quasiconformal maps induced by different types of discrete conformal structures converge to a $1$-quasiconformal map, i.e., conformal map.

For general discrete conformal structures, we propose a similar infinite rigidity theorem.
\begin{conjecture}\label{Conjecture Euclidean}
Let $(\mathbb{C}, \mathcal{T}, \varepsilon, \eta)$ be a triangulated plane with a regular weight pair $(\varepsilon, \eta)$.
If $l$ and $\bar{l}$ are two flat weighted Delaunay PL metrics on $(\mathbb{C}, \mathcal{T}, \varepsilon, \eta)$, i.e., $\bar l=w^e*l$, 
then $w^e \equiv C$ for some constant $C$.
\end{conjecture}

Conjecture \ref{Conjecture Euclidean} represents a discrete analogue of the characterization of conformal automorphisms of the complex plane: every conformal map from the plane to itself is a similarity.

\begin{conjecture}
Let $(\mathbb{D}, \mathcal{T}, \varepsilon, \eta)$ be a triangulated plane with a regular weight pair $(\varepsilon, \eta)$.
If $l$ and $\bar{l}$ are two flat weighted Delaunay PH metrics on $(\mathbb{D}, \mathcal{T}, \varepsilon, \eta)$, i.e., $\bar{l}=w^h*l$,
then $w^h \equiv 0$, i.e., they are isometric.
\end{conjecture}


\subsection{Organization of the paper}

In Section \ref{Sec: DCS-WDT}, we review some fundamental properties of the Euclidean discrete conformal structures and weighted Delaunay triangulations.
Subsequently, we extend the definition of weighted Delaunay triangulations to encompass generalized triangles and establish several useful lemmas.
In Section \ref{Sec: MP-E}, we prove the Euclidean discrete maximum principle, i.e., Theorem \ref{Thm: MP-E}.
In Section \ref{Sec: MP-H}, we derive the hyperbolic discrete maximum principle based on the relationships between hyperbolic and Euclidean discrete conformal structures.
As an application of the hyperbolic discrete maximal principle, we further prove the discrete Schwarz-Ahlfors lemma in this section.
In Section \ref{Sec: IR}, we prove the infinite rigidity theorem for small Delaunay triangulations of the hyperbolic plane.
\\
\\
\textbf{Acknowledgment}\\[8pt]
The authors thank Dr. Tianqi Wu for the insight on Theorem \ref{Thm: IR-H}.
The research of X. Xu is supported by National Natural Science Foundation of China
under grant no. 12471057.

\section{Euclidean discrete conformal structures and weighted Delaunay triangulations}
\label{Sec: DCS-WDT}

In this section, we first give some basic properties of Euclidean discrete conformal structures, then we generalize the definition of weighted Delaunay for nondegenerate triangles generated by Euclidean discrete conformal structures to generalized triangles.
Given that this section exclusively addresses the Euclidean case and does not consider the hyperbolic case, the term ``Euclidean" is omitted for simplicity.

\subsection{Basic properties of Euclidean discrete conformal structures}

To be consist with the notations used in the hyperbolic case, we define 
\begin{equation*}
u=f = \ln r.
\end{equation*}
By combining $\varepsilon \in \{0,1\}$ with (\ref{Eq: SC1}), we obtain

\begin{equation*}
\varepsilon_ie^{2f_i}+\varepsilon_je^{2f_j}
+2\eta_{ij}e^{f_i+f_j}
\geq 2(\varepsilon_i\varepsilon_j+\eta_{ij})e^{f_i+f_j}>0.
\end{equation*}
Thus $l_{ij}$ in (\ref{Eq: length-E}) is well-defined.
The lengths $l_{12},l_{23},l_{13}$ forms a non-degenerate triangle $\triangle v_1v_2v_3$ if and only if they satisfy the following strict triangle inequality
\begin{equation*}
l_{ij} < l_{jk} + l_{ki}, \quad \{i,j,k\} = \{1,2,3\}.
\end{equation*}
In this case, we say the labels $f,r,u$ are admissible on $\triangle v_1v_2v_3$. 
Let $\Omega_{123}(\eta)$ denote the space of admissible labels for $\triangle v_1v_2v_3$, parametrized by $f=u$ or $r$, if it  does not cause confusion in the context. 
If all triangles in $(S, \mathcal{T}, \varepsilon, \eta)$ are non-degenerate, we say that $f,r,$ or $u$ are admissible labels on $(S, \mathcal{T}, \varepsilon, \eta)$.

If the lengths $l_{ij},l_{ik},l_{jk}$ satisfy the triangle inequality
\begin{equation*}
l_{ij} \leq l_{jk} + l_{ki}, \quad \{i,j,k\} = \{1,2,3\},
\end{equation*}
then we refer to $\triangle v_1v_2v_3$ as a generalized triangle, and $f,r,u$ as generalized labels. 
If all the triangles in $(S, \mathcal{T}, \varepsilon, \eta)$ are generalized triangles, we say that $f,r,u$ are generalized labels on $(S, \mathcal{T}, \varepsilon, \eta)$. 
In particular, if $l_{ij}=l_{ik}+l_{kj}$, we refer to the generalized triangle $\triangle v_1v_2v_3$ as a degenerate triangle. 
In this case, the inner angle at $v_k$ is defined to be $\pi$, and we call $v_k$ the flat vertex of $\triangle v_1v_2v_3$.

\begin{lemma}[\cite{Xu 1}]\label{Lem: basic 1}
Suppose $(S,\mathcal{T},\varepsilon,\eta)$ is a weighted triangulated surface with $\varepsilon: V\rightarrow \{0,1\}$ and $\eta: E\rightarrow \mathbb{R}_{\geq 0}$ satisfying the structure condition  (\ref{Eq: SC1}).
Let $\triangle v_1v_2v_3$ be a generalized triangle with edge lengths defined by (\ref{Eq: length-E}) or (\ref{Eq: length-E2}) on $(S, \mathcal{T}, \varepsilon, \eta)$.
\begin{description}
\item[(a)]
The triangle $\triangle v_1v_2v_3$ is non-degenerate if and only if $Q>0$, where
\begin{equation}\label{Eq: Q}
Q=\kappa_1^2(\varepsilon_2\varepsilon_3-\eta^2_{23})
+\kappa_2^2(\varepsilon_1\varepsilon_3-\eta^2_{13})
+\kappa_3^2(\varepsilon_1\varepsilon_2-\eta^2_{12}) +2\kappa_1\kappa_2\gamma_{3}
+2\kappa_1\kappa_3\gamma_{2}
+2\kappa_2\kappa_3\gamma_{1}
\end{equation}
with $\gamma_{i}:=\varepsilon_i\eta_{jk}+\eta_{ij}\eta_{ik}\geq 0$ and $\kappa_i:=r_i^{-1}$ for $\{i,j,k\}=\{1,2,3\}$.
As a result, $\triangle v_1v_2v_3$ is degenerate if and only if $Q=0$.

\item[(b)]
Let $\theta_i$ be the inner angle at vertex $v_i$.
If $\triangle v_1v_2v_3$ is non-degenerate,
then
\begin{equation*}
\frac{\partial \theta_i}{\partial u_j}
=\frac{\partial \theta_j}{\partial u_i}
=\frac{h_{ij,k}}{l_{ij}},\
\frac{\partial \theta_i}{\partial u_i}
=-\frac{\partial \theta_i}{\partial u_j}
-\frac{\partial \theta_i}{\partial u_k}
<0,
\end{equation*}
where 
\begin{equation}\label{Eq: hij,k}
h_{ij,k}
=\frac{r_1^2r_2^2r_3^2}{A_{123}l_{ij}}
[\kappa_k^2(\varepsilon_i\varepsilon_j-\eta_{ij}^2)
+\kappa_j\kappa_k\gamma_{i}+\kappa_i\kappa_k\gamma_{j}]				=\frac{r_1^2r_2^2r_3^2}{A_{123}l_{ij}}\kappa_kh_k,
\end{equation}
with $A_{123}=l_{12}l_{23}\sin \theta_2$, and
\begin{equation}\label{Eq: h_i}						h_k=\kappa_k(\varepsilon_i\varepsilon_j-\eta_{ij}^2)
+\kappa_j\gamma_i+\kappa_i\gamma_j.
\end{equation}
Furthermore, the Jacobian $\Lambda_{ijk}:=\frac{\partial (\theta_1, \theta_2, \theta_3)}{\partial (u_1, u_2, u_3)}$ is symmetric and negative semi-definite with rank 2 and kernel $\{c(1,1,1)^T|c\in \mathbb{R} \}$ on $\Omega_{123}(\eta)$.

\item[(c)]
If $\triangle v_1v_2v_3$ induced by $(r_1, r_2, r_3)$ via (\ref{Eq: length-E2}) is degenerate, 
then one of $h_1$, $h_2$, and $h_3$ is negative, while the other two are positive. 
Specially, if $v_3$ is the flat vertex of $\triangle v_1v_2v_3$, 
then $h_1>0$, $h_2>0$, and $h_3<0$. 
This implies that $\eta^2_{12}>\varepsilon_1\varepsilon_2$ due to $h_3<0$, equivalently $\eta_{12}>\varepsilon_1\varepsilon_2$ by (\ref{Eq: SC1}). Furthermore, if $v_3$ is the flat vertex of $\triangle v_1v_2v_3$ induced by $(r_1, r_2, r_3)$, as $(\tilde{r}_1, \tilde{r}_2, \tilde{r}_3) \in \Omega_{123}(\eta)$ tends to $(r_1, r_2, r_3) \in \partial \Omega_{123}(\eta)$, we have
\begin{equation*}
h_{12,3} \rightarrow -\infty, \quad
h_{13,2} \rightarrow +\infty, \quad
h_{23,1} \rightarrow +\infty.
\end{equation*}

\item[(d)]
Define
\begin{equation}\label{Eq: F1}
\Delta_{ijk}
=\varepsilon_i\eta_{jk}^2+\varepsilon_j\eta_{ik}^2
+\varepsilon_k\eta_{ij}^2+2\eta_{ij}\eta_{ik}\eta_{jk}
-\varepsilon_i\varepsilon_j\varepsilon_k.
\end{equation}
If at least one of the inequalities $\eta_{12}>\varepsilon_1\varepsilon_2$, $\eta_{13}>\varepsilon_1\varepsilon_3$, or $\eta_{23}>\varepsilon_2\varepsilon_3$ holds, then $\Delta_{123}>0$. Consequently, if $\Delta_{123}\leq 0$, then $\eta_{12}\leq\varepsilon_1\varepsilon_2$, $\eta_{13}\leq\varepsilon_1\varepsilon_3$, and $\eta_{23}\leq\varepsilon_2\varepsilon_3$. This implies that $\varepsilon_1=\varepsilon_2=\varepsilon_3=1$ by (\ref{Eq: SC1}). 
Furthermore, in the case $\Delta_{123}\leq 0$,
it follows from (\ref{Eq: Q}) that $Q(r)>0$ for any $r\in\mathbb{R}_{>0}^3$.
\end{description}
\end{lemma}

\begin{lemma}\label{Lem: basic 2}
Suppose $(S,\mathcal{T},\varepsilon,\eta)$ is a weighted triangulated surface with $\varepsilon: V\rightarrow \{0,1\}$ and $\eta: E\rightarrow \mathbb{R}_{\geq 0}$ satisfying the structure conditions (\ref{Eq: SC1}).
Let $\triangle v_1v_2v_3$ be a generalized triangle induced by $(r_1,r_2,r_3)\in\mathbb{R}^3_{>0}$ via (\ref{Eq: length-E2}) on $(S,\mathcal{T},\varepsilon,\eta)$.
\begin{description}
\item[(a)]
The admissible space $\Omega_{123}(\eta)$ of admissible labels $(r_1, r_2, r_3) \in \mathbb{R}^3_{>0}$ is a non-empty, simply connected open subset of $\mathbb{R}^3_{>0}$ with analytic boundary components.
Furthermore,
\begin{equation*}
\Omega_{123}(\eta)=\mathbb{R}^3_{>0}\setminus \bigsqcup_{\alpha \in \Lambda} V_\alpha,
\end{equation*}
where $\Lambda=\{i \mid A_i=\eta^2_{jk}-\varepsilon_j \varepsilon_k>0, \ \{i,j,k\}=\{1,2,3\}\}$,
and $\sqcup_{\alpha \in \Lambda} V_\alpha$ is a disjoint union of
\begin{equation*}
V_i = \left\{(r_1, r_2, r_3) \in \mathbb{R}^3_{>0} \mid \kappa_i \geq \frac{-B_i + \sqrt{\Delta_i}}{2A_i}\right\}
\end{equation*}
with
\begin{equation*}
\begin{aligned}				B_i&=-2(\kappa_j\gamma_{k}+\kappa_{k}\gamma_j)\leq0,\\
\Delta_i&=4(\varepsilon_k\kappa_j^2
+\varepsilon_j\kappa_{k}^2+2\kappa_j\kappa_{k}\eta_{jk})
\cdot\Delta_{ijk}\\
&=4(\varepsilon_k\kappa_j^2
+\varepsilon_j\kappa_{k}^2+2\kappa_j\kappa_{k}\eta_{jk})\cdot (\varepsilon_i\eta_{jk}^2+\varepsilon_j\eta_{ik}^2
+\varepsilon_k\eta_{ij}^2+2\eta_{ij}\eta_{ik}\eta_{jk}
-\varepsilon_i\varepsilon_j\varepsilon_k).
\end{aligned}
\end{equation*}

\item[(b)]
The inner angles of $\triangle v_1v_2v_3$ could be uniquely continuously extended by constants as follows
\begin{equation*}
\begin{aligned}
\tilde{\theta}_i(r_1, r_2, r_3) = \left\{
\begin{array}{ll}
\theta_i, & \text{if } (r_1, r_2, r_3) \in \Omega_{123}(\eta), \\
\pi, & \text{if } (r_1, r_2, r_3) \in V_i, \\
0, & \text{otherwise}.
\end{array}
\right.
\end{aligned}
\end{equation*}

\item[(c)]
If $\triangle v_1v_2v_3$ induced by $(r_1, r_2, r_3)$ is degenerate with $v_i$ as the flat vertex,
then $(r_1, r_2, r_3) \in \partial V_i$, i.e.,
$(r_1, r_2, r_3)$ satisfies
\begin{equation}\label{Eq: F11}
\kappa_i = \frac{-B_i + \sqrt{\Delta_i}}{2A_i}.
\end{equation}
Furthermore, as a function of $\kappa_j$ and $\kappa_k$, $\kappa_i := f(\kappa_j, \kappa_k)$ is strictly increasing with respect to both $\kappa_j$ and $\kappa_k$.
\end{description}
\end{lemma}
\proof
Parts (a) and (b) have been proved in \cite{Xu 1}, it remains to prove part (c). 
Part (a) shows that $\partial\Omega_{123}(\eta) = \partial V_1 \sqcup \partial V_2 \sqcup \partial V_3$. 
The assumption of part (c) implies that the inner angle at $v_i$ is $\pi$. 
By part (b), we have $(r_1,r_2,r_3)\in\partial V_i$, 
which yields (\ref{Eq: F11}). 
Note that (\ref{Eq: F11}) can be explicitly rewritten as
\begin{gather*}
f(\kappa_j, \kappa_k) := \kappa_i
= \frac{1}{\eta^2_{jk} - \varepsilon_j \varepsilon_k}
\left[(\kappa_j \gamma_{k} + \kappa_{k} \gamma_j) +
\sqrt{(\varepsilon_k \kappa_j^2
+ \varepsilon_j \kappa_{k}^2 + 2\kappa_j \kappa_{k} \eta_{jk})
\cdot \Delta_{ijk}}\right].
\end{gather*}
Combining the assumption with Lemma \ref{Lem: basic 1} (c), we obtain $\eta_{jk}>\varepsilon_j\varepsilon_k$. 
From Lemma \ref{Lem: basic 1} (d), it follows that $\Delta_{ijk}>0$. 
Since $\varepsilon_k\kappa_j^2 +\varepsilon_j\kappa_{k}^2+2\kappa_j\kappa_{k}\eta_{jk} \geq2(\varepsilon_j\varepsilon_k+\eta_{jk})\kappa_j\kappa_{k} >0$ by (\ref{Eq: SC1}), it follows that $f(\kappa_j,\kappa_k)>0$ for any $(\kappa_j,\kappa_k)\in\mathbb{R}^2_{>0}$. 
It is straightforward to check that $f(\kappa_j,\kappa_k)$ is strictly increasing with respect to both $\kappa_j$ and $\kappa_k$ by the structure condition  (\ref{Eq: SC1}).
\qed

\begin{remark}
Note that Lemma \ref{Lem: basic 1} and Lemma \ref{Lem: basic 2} hold for $\varepsilon: V\rightarrow \{0,1\}$ and $\eta: E\rightarrow \mathbb{R}$ satisfying the structure conditions (\ref{Eq: SC1}) and
\begin{equation}\label{Eq: SC2}
\gamma_{i}\geq 0, \quad \gamma_j\geq 0, \quad \gamma_k\geq 0
\end{equation}
for any triangle $\triangle v_iv_jv_k\in F$. 
In particular, if $\eta: E\rightarrow \mathbb{R}_{\geq 0}$, then (\ref{Eq: SC2}) is automatically satisfied.
\end{remark}

\begin{lemma}\label{Lem: interval}
Suppose $(S,\mathcal{T},\varepsilon,\eta)$ is a weighted triangulated surface with $\varepsilon: V\rightarrow \{0,1\}$ and $\eta: E\rightarrow \mathbb{R}_{\geq 0}$ satisfying the structure condition (\ref{Eq: SC1}).
Let $\triangle v_1v_2v_3$ be a generalized triangle induced by $(r_1,r_2,r_3)\in\mathbb{R}^3_{>0}$ via (\ref{Eq: length-E2}) on $(S,\mathcal{T},\varepsilon,\eta)$.
\begin{description}
\item[(a)]
For any fixed $r_i, r_j\in \mathbb{R}^2_{>0}$,
the set of $r_k\in (0, +\infty)$ such that the label $(r_1,r_2,r_3)$ is admissible is an open interval.
As a result, if $(r_i, r_j,\hat{r}_k)$ and $(r_i, r_j,\bar{r}_k)$ are two generalized labels on $\triangle v_1v_2v_3$ with $\hat{r}_k<\bar{r}_k$,
then for any $r_k\in (\hat{r}_k,\bar{r}_k)$,
the triangle $\triangle v_1v_2v_3$ induced by $(r_i, r_j, r_k)$ is non-degenerate.

\item[(b)]
If $\triangle v_1v_2v_3$ induced by $(r_1,r_2,r_3)$ is degenerate, and $v_3$ is the flat vertex of $\triangle v_1v_2v_3$,
then there exists $\epsilon>0$ such that $(r_1,r_2,r_3+t)\in \Omega_{123}(\eta)$, and
\begin{equation*}
\frac{\partial h_{12,3}}{\partial r_3}(r_1, r_2, r_3+t)>0
\end{equation*}
for $t\in (0, \epsilon)$.
\end{description}
\end{lemma}
\proof 
\textbf{(a)}
Without loss of generality, let $\{i,j\}=\{2,3\}$ and $k=1$. 
The quantity $Q$ in (\ref{Eq: Q}) can be expressed as a quadratic function of $\kappa_1$
\begin{equation*}
g(\kappa_1)
=(\varepsilon_2\varepsilon_3-\eta^2_{23})\kappa_1^2
+2(\kappa_2\gamma_{3}+\kappa_3\gamma_{2})\kappa_1
+\kappa_2^2(\varepsilon_1\varepsilon_3-\eta^2_{13})
+\kappa_3^2(\varepsilon_1\varepsilon_2-\eta^2_{12})
+2\kappa_2\kappa_3\gamma_{1}.
\end{equation*}
By Lemma \ref{Lem: basic 1} (a), it suffices to show that the solution set of $g(\kappa_1)>0$ for $\kappa_1\in (0, +\infty)$ is an open interval. It is necessary to consider the following three cases $\eta^2_{23}=\varepsilon_2\varepsilon_3$, $\eta^2_{23}<\varepsilon_2\varepsilon_3$ and $\eta^2_{23}>\varepsilon_2\varepsilon_3$.

Assume that $\eta^2_{23}=\varepsilon_2\varepsilon_3$, i.e., $(\eta_{23}-\varepsilon_2\varepsilon_3)(\eta_{23}+\varepsilon_2\varepsilon_3)=0$. 
From (\ref{Eq: SC1}), it follows that $\eta_{23}=\varepsilon_2=\varepsilon_3=1$. Consequently, $g(\kappa_1)>0$ is equivalent to
\begin{equation}\label{Eq: F2}
g(\kappa_1)=
2(\kappa_2\gamma_{3}+\kappa_3\gamma_{2})\kappa_1
+\kappa_2^2(\varepsilon_1-\eta^2_{13})
+\kappa_3^2(\varepsilon_1-\eta^2_{12})
+2\kappa_2\kappa_3\gamma_{1}>0.
\end{equation}
If $\gamma_{2}=\gamma_{3}=0$, then $\gamma_{2}+\gamma_{3}=(\eta_{12}+\eta_{13})(\eta_{23}+1)=0$, and thus $\eta_{12}=\eta_{13}=0$ by $\eta\in \mathbb{R}^E_{\geq 0}$. From (\ref{Eq: SC1}), it follows that $\varepsilon_1=1$. Consequently, (\ref{Eq: F2}) reduces to
\begin{equation*}
g(\kappa_1)=
\kappa_2^2+\kappa_3^2+2\kappa_2\kappa_3
=(\kappa_2+\kappa_3)^2>0.
\end{equation*}
Therefore, the solution set of $g(\kappa_1)>0$ is $\mathbb{R}_{>0}$ in this case. If at least one of $\gamma_{2}$ and $\gamma_{3}$ is positive, then $\kappa_2\gamma_{3}+\kappa_3\gamma_{2}>0$. Therefore, the solution of (\ref{Eq: F2}) is
\begin{equation*}
\kappa_1>-\frac{\kappa_2^2(\varepsilon_1-\eta^2_{13})
+\kappa_3^2(\varepsilon_1-\eta^2_{12})
+2\kappa_2\kappa_3\gamma_{1}}
{2(\kappa_2\gamma_{3}+\kappa_3\gamma_{2})}.
\end{equation*}
This implies that the solution set of $g(\kappa_1)>0$ for $\kappa_1\in (0, +\infty)$ is an open interval in this case.

Assume that $\eta^2_{23}<\varepsilon_2\varepsilon_3$. Then $g(\kappa_1)$ is a quadratic function opening upward with respect to $\kappa_1$. 
Note that the axis of symmetry of the quadratic function $g(\kappa_1)$ is given by $-\frac{\kappa_2\gamma_{3}+\kappa_3\gamma_{2}}{\varepsilon_2\varepsilon_3-\eta^2_{23}}\leq0$. Therefore, the solution set of $g(\kappa_1)>0$ for $\kappa_1\in (0, +\infty)$ is an open interval in this case.

Assume that $\eta^2_{23}>\varepsilon_2\varepsilon_3$. 
Then $g(\kappa_1)$ is a quadratic function opening downward with respect to $\kappa_1$.
Note that the axis of symmetry of the quadratic function $g(\kappa_1)$ is given by $-\frac{\kappa_2\gamma_{3}+\kappa_3\gamma_{2}}{\varepsilon_2\varepsilon_3-\eta^2_{23}}\geq0$, and the discriminant of $g(\kappa_1)$ is
\begin{equation*}
\Delta=4(\varepsilon_2\kappa_3^2
+\varepsilon_3\kappa_2^2+2\kappa_2\kappa_3\eta_{23})
\cdot\Delta_{123}.
\end{equation*}
Since $\eta^2_{23}>\varepsilon_2\varepsilon_3$, 
it follows from (\ref{Eq: SC1}) that $\eta_{23}>\varepsilon_2\varepsilon_3$. 
From Lemma \ref{Lem: basic 1} (d), we have $\Delta_{123}>0$. 
Furthermore, $\varepsilon_2\kappa_3^2+\varepsilon_3\kappa_2^2+2\kappa_2\kappa_3\eta_{23}\geq 2(\varepsilon_2\varepsilon_3+\eta_{23})\kappa_2\kappa_3>0$ by (\ref{Eq: SC1}). 
Therefore, $\Delta>0$. 
Consequently, the solution set of $g(\kappa_1)>0$ for $\kappa_1\in (0, +\infty)$ is an open interval in this case.

\textbf{(b)} 
From Lemma \ref{Lem: basic 1} (a), it follows that $\triangle v_1v_2v_3$ is degenerate if and only if $Q=0$. 
Direct calculations yield 
\begin{equation*}
\frac{\partial Q}{\partial r_3}=\frac{\partial Q}{\partial \kappa_3}\frac{\partial \kappa_3}{\partial r_3}=-2\kappa^2_3 h_3,
\end{equation*}
where $h_3$ is defined by (\ref{Eq: h_i}).
Since $v_3$ is the flat vertex of $\triangle v_1v_2v_3$, it follows from Lemma \ref{Lem: basic 1} (c) that $h_3<0$ at $(r_1,r_2,r_3)$. 
This implies that $\frac{\partial Q}{\partial r_3}>0$ in a neighborhood of $(r_1, r_2, r_3)$. 
Therefore, for sufficiently small $t>0$, $Q(r_1,r_2,r_3+t)>0$ and $(r_1,r_2,r_3+t)$ induces a non-degenerate triangle. 

Direct calculations yield
\begin{equation}\label{Eq: F10}
A_{123}^2=r_1^2r_2^2r_3^2Q.
\end{equation}
For a detailed derivation of (\ref{Eq: F10}), please refer to the arguments above Lemma 2.1 in \cite{Xu 1}. Consequently,
\begin{equation*}
\frac{\partial A_{123}}{\partial r_{3}}
=\frac{A_{123}}{r_{3}}-\frac{r_{1}^{2} r_{2}^{2}}{A_{123}} \cdot h_{3}.
\end{equation*}
By differentiating (\ref{Eq: hij,k}) with respect to $r_3$ and combining the result with (\ref{Eq: F10}), we derive
\begin{equation}\label{Eq: F3}
\frac{\partial h_{12,3}}{\partial r_{3}}
=\frac{r_1^4r_2^4r_3}{A_{123}^3l_{12}}\Delta_{123}
(\varepsilon_2\kappa_1^2+\varepsilon_1\kappa_2^2
+2\kappa_1\kappa_2\eta_{12}).
\end{equation}
In fact,
\begin{equation*}
\begin{aligned}
\frac{\partial h_{12,3}}{\partial r_{3}} & =\frac{r_{1}^{2} r_{2}^{2}}{l_{12}} \frac{\partial}{\partial r_{3}}\left(\frac{r_{3} h_{3}}{A_{123}}\right)\\
&=\frac{r_{1}^{2} r_{2}^{2}}{l_{12}} \cdot \frac{1}{A_{123}^{2}}\left[\left(h_{3}+r_{3} \cdot \frac{\partial h_{3}}{\partial r_{3}}\right) \cdot A_{123}-r_{3} h_{3} \cdot \frac{\partial A_{123}}{\partial r_{3}}\right] \\
& =\frac{r_{1}^{2} r_{2}^{2}}{l_{12} \cdot A_{123}^{2}}\left[h_{3} A_{123}+r_{3} A_{123} \cdot\left(-\frac{1}{r_{3}^{2}}\right)\left(\varepsilon_{1} \varepsilon_{2}-\eta_{12}^{2}\right)-r_{3} h_{3} \cdot\left(\frac{A_{123}}{r_{3}}-\frac{r_{1}^{2} r_{2}^{2}}{A_{123}} \cdot h_{3}\right)\right] \\
& =\frac{r_{1}^{2} r_{2}^{2}}{l_{12} A_{123}^{2}}\left[-A_{123} \kappa_{3}\left(\varepsilon_{1} \varepsilon_{2}-\eta_{12}^{2}\right)+\frac{r_{1}^{2} r_{2}^{2} r_{3}}{A_{123}} \cdot h_{3}^{2}\right] \\
& =\frac{r_{1}^4 r_{2}^4r_3}{l_{12} A_{123}^3}\left[-\left(\varepsilon_{1} \varepsilon_{2}-\eta_{12}^{2}\right)Q+h_{3}^{2}\right] \\
&=\frac{r_1^4r_2^4r_3}{A_{123}^3l_{12}}\Delta_{123}
(\varepsilon_2\kappa_1^2+\varepsilon_1\kappa_2^2
+2\kappa_1\kappa_2\eta_{12}).
\end{aligned}
\end{equation*}
In the fifth line, substituting the explicit expressions of $Q$ in (\ref{Eq: Q}) and $h_3$ in (\ref{Eq: h_i}) and then simplifying yields the last line.

Since $\triangle v_1v_2v_3$ induced by $(r_1,r_2,r_3)$ is degenerate with $v_3$ as the flat vertex, 
it follows from Lemma \ref{Lem: basic 1} (c) that $\eta_{12}>\varepsilon_1\varepsilon_2$. 
From Lemma \ref{Lem: basic 1} (d), we have $\Delta_{123}>0$. 
Furthermore, $\varepsilon_2\kappa_1^2+\varepsilon_1\kappa_2^2 +2\kappa_1\kappa_2\eta_{12}>0$. 
Note that $h_{12,3}$ is only defined for admissible labels. 
Therefore, these inequalities hold in a neighborhood of $(r_1, r_2, r_3)$. 
Consequently, there exists $\epsilon>0$ such that $\frac{\partial h_{12,3}}{\partial r_3}(r_1, r_2, r_3+t)>0$ for $t\in (0, \epsilon)$.
\qed

\subsection{Weighted Delaunay triangulations and some useful lemmas}
\label{Sec: weighted Delaunay}

Let $\triangle v_1v_2v_3$ be a non-degenerate triangle induced by $(r_1,r_2,r_3)\in\mathbb{R}^3_{>0}$. 
The projections of the center $c_{123}$ of the face-circle $C_{123}$ onto the lines determined by $v_1v_2$, $v_1v_3$ and $v_2v_3$ yield the geometric centers of these edges, denoted by $c_{12}$, $c_{13}$ and $c_{23}$, respectively. 
The signed distance from $c_{12}$ to $v_1$, denoted by $d_{12}$, is positive if $c_{12}$ lies on the same side as $v_2$ along the line defined by $v_1v_2$, negative otherwise, and zero if $c_{12}$ coincides with $v_1$. 
The signed distance $d_{21}$ is defined in a similar way. 
Glickenstein \cite{Glickenstein JDG} obtained the following identities
\begin{equation}\label{Eq: d}
d_{ij}=\frac{\varepsilon_ir_i^2+r_ir_j\eta_{ij}}{l_{ij}}, \quad h_{ij,k}=\frac{d_{ik}-d_{ij}\cos \theta_i}{\sin \theta_i}.
\end{equation}
Combining (\ref{Eq: SC1}) with (\ref{Eq: d}) yields $d_{ij}>0$. 
Note that $d_{ij}+d_{ji}=l_{ij}$, and $d_{ij}\neq d_{ji}$ in general.

\begin{definition}\label{Def: theta}
Suppose $(S,\mathcal{T},\varepsilon,\eta)$ is a weighted triangulated surface with $\varepsilon: V\rightarrow \{0,1\}$ and $\eta: E\rightarrow \mathbb{R}_{\geq 0}$ satisfying the structure condition (\ref{Eq: SC1}).
Let $\triangle v_1v_2v_3\in \mathcal{T}$ be a generalized triangle induced by a generalized label.
If $\triangle v_1v_2v_3$ is non-degenerate, we define
\begin{equation}\label{Eq: F6}
\theta_{ij,k}=\arctan\frac{h_{ij,k}}{d_{ij}}.
\end{equation}
If $\triangle v_1v_2v_3$ is degenerate, we define 
\begin{equation}\label{Eq: F7}
\theta_{ij,k}=
\begin{cases} 
+\frac{\pi}{2}, & \text{if $v_i$ or $v_{j}$ is the flat vertex}, \\
-\frac{\pi}{2}, & \text{if $v_k$ is the flat vertex}.
\end{cases}
\end{equation}
\end{definition}

\begin{figure}[!ht]
\centering
\includegraphics[scale=1]{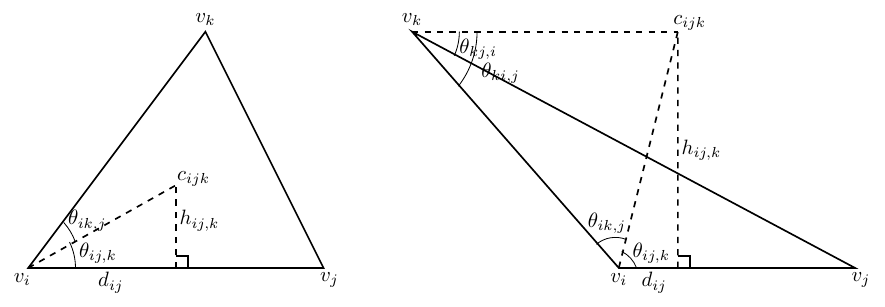}
\caption{The schematic illustration of the angle $\theta_{ij,k}$.}
\label{figure5}
\end{figure}

For a non-degenerate triangle $\triangle v_1v_2v_3$, $\theta_{ij,k}$ is exactly the signed angle $\angle v_jv_ic_{123}$, which is positive if $c_{123}$ is on the same side of the line determined by $v_iv_j$ as the triangle $\triangle v_1v_2v_3$,
negative otherwise, and zero if $c_{123}$ lies in the line determined by $v_iv_j$.
It is straight forward to check that $\theta_{ij,k}$ is a continuous function of $(r_1,r_2,r_3)\in \Omega_{123}$.
Moreover, $\theta_{ij,k}+\theta_{ik,j}=\theta_i$.
Please refer to Figure \ref{figure5}.

\begin{lemma}\label{Lem: theta continuous}
Suppose $(S,\mathcal{T},\varepsilon,\eta)$ is a weighted triangulated surface with $\varepsilon: V\rightarrow \{0,1\}$ and $\eta: E\rightarrow \mathbb{R}_{\geq 0}$ satisfying the structure condition (\ref{Eq: SC1}).
Let $\triangle v_1v_2v_3$ be a generalized triangle induced by a generalized label $(r_1, r_2, r_3)$ on $(S,\mathcal{T},\varepsilon,\eta)$.
Then $\theta_{ij,k}(r_1, r_2, r_3)$ is a continuous function defined on  $\overline{\Omega}_{123}(\eta)$, and satisfies
\begin{equation}\label{Eq: F4}
\theta_{ij,k}+\theta_{ik,j}=\theta_i.
\end{equation}
\end{lemma}
\proof
The proof is nearly identical to that of Lemma 2.8 in \cite{LXZ}, and we provide a proof here for completeness. 
It suffices to show that as $(r_1,r_2,r_3)\in \Omega_{123}(\eta)$ tends to $(\bar{r}_1,\bar{r}_2,\bar{r}_3)\in \partial\Omega_{123}(\eta)$,  $\theta_{ij,k}(r_1,r_2,r_3)\rightarrow \theta_{ij,k}(\bar{r}_1,\bar{r}_2,\bar{r}_3)$. 
Assume that $v_k$ is the flat vertex of the degenerate triangle $\triangle v_1v_2v_3$ induced by $(\bar{r}_1,\bar{r}_2,\bar{r}_3)$. 
As $(r_1,r_2,r_3)\rightarrow(\bar{r}_1,\bar{r}_2,\bar{r}_3)$, it follows from Lemma \ref{Lem: basic 1} (c) that $h_{ij,k}(r_1,r_2,r_3)\rightarrow -\infty$. 
This implies that $\theta_{ij,k}(r_1,r_2,r_3)=\arctan\frac{h_{ij,k}}{d_{ij}}\rightarrow -\frac{\pi}{2}=\theta_{ij,k}(\bar{r}_1,\bar{r}_2,\bar{r}_3)$ by Definition \ref{Def: theta}. 
Assume that $v_i$ is the flat vertex of the degenerate triangle $\triangle v_1v_2v_3$ induced by $(\bar{r}_1,\bar{r}_2,\bar{r}_3)$. 
Similarly, as $(r_1,r_2,r_3)\rightarrow(\bar{r}_1,\bar{r}_2,\bar{r}_3)$, it follows that $h_{ij,k}(r_1,r_2,r_3)\rightarrow +\infty$. 
This implies that $\theta_{ij,k}(r_1,r_2,r_3)\rightarrow \frac{\pi}{2}=\theta_{ij,k}(\bar{r}_1,\bar{r}_2,\bar{r}_3)$ by Definition \ref{Def: theta}. 
The same argument holds for the case that $v_j$ is the flat vertex.
\qed

\begin{corollary}\label{Cor: weighted Delaunay 2}
Suppose $(S,\mathcal{T},\varepsilon,\eta)$ is a weighted triangulated surface with $\varepsilon: V\rightarrow \{0,1\}$ and $\eta: E\rightarrow \mathbb{R}_{\geq 0}$ satisfying the structure condition (\ref{Eq: SC1}).
Let $r\in \mathbb{R}^V_{>0}$ be an admissible label on $(S,\mathcal{T},\varepsilon,\eta)$.
For any two adjacent non-degenerate triangles $\triangle v_iv_jv_k$ and $\triangle v_iv_jv_l$ sharing the common edge $v_iv_j$,
the edge $v_iv_j$ is weighted Delaunay in $r$ if and only if
\begin{equation*}
\theta_{ij,k}+\theta_{ij,l}\ge0.
\end{equation*}
\end{corollary}
\proof
From Definition \ref{Def: theta}, it follows that
\begin{equation*}
\frac{h_{ij,k}+h_{ij,l}}{d_{ij}}
=\tan\theta_{ij,k}+\tan\theta_{ij,l}
=\frac{\sin(\theta_{ij,k}+\theta_{ij,l})}
{\cos\theta_{ij,k}\cos\theta_{ij,l}}.
\end{equation*}
Note that $d_{ij}>0$ and $\theta_{ij,k}, \theta_{ij,l}\in (-\frac{\pi}{2},\frac{\pi}{2})$ for non-degenerate triangles $\triangle v_iv_jv_k$ and $\triangle v_iv_jv_l$,
we have $h_{ij,k}+h_{ij,l}\ge0$ is equivalent to $\theta_{ij,k}+\theta_{ij,l}\ge0$.
\qed

Lemma \ref{Lem: theta continuous} and Corollary \ref{Cor: weighted Delaunay 2} motivate the following definition of weighted Delaunay triangulation for
generalized discrete conformal structures on surfaces.

\begin{definition}\label{Def: weighted Delaunay-G}
Suppose $(S,\mathcal{T},\varepsilon,\eta)$ is a weighted triangulated surface with $\varepsilon: V\rightarrow \{0,1\}$ and $\eta: E\rightarrow \mathbb{R}_{\geq 0}$ satisfying the structure condition (\ref{Eq: SC1}).
Let $r\in \mathbb{R}^V_{>0}$ be a generalized label on $(S,\mathcal{T},\varepsilon,\eta)$.
For any two adjacent generalized triangles $\triangle v_iv_jv_k$ and $\triangle v_iv_jv_l$ sharing the common edge $v_iv_j$,
the edge $v_iv_j$ is weighted Delaunay in $r$ if
\begin{equation}\label{Eq: F8}
\theta_{ij,k}+\theta_{ij,l}\ge0.
\end{equation}
The triangulation $\mathcal{T}$ is weighted Delaunay in $r$ if (\ref{Eq: F8}) holds for every interior edge in $\mathcal{T}$.
\end{definition}

The following two lemmas play important roles in the proof of Lemma \ref{Lem: key}.

\begin{lemma}\label{Lem: monotonicity}
Suppose $(S,\mathcal{T},\varepsilon,\eta)$ is a weighted triangulated surface with $\varepsilon: V\rightarrow \{0,1\}$ and $\eta: E\rightarrow \mathbb{R}_{\geq 0}$ satisfying the structure condition (\ref{Eq: SC1}).
Let $(r_1, r_2, \hat{r}_3)$ and $(r_1, r_2, \bar{r}_3)$ be two generalized labels on $\triangle v_1v_2v_3$ with $\hat{r}_3<\bar{r}_3$.
If $\Delta_{123}>0$, then for fixed $r_1$ and $r_2$, $\theta_{12,3}$ is strictly increasing with respect to $r_3\in [\hat{r}_3,\bar{r}_3]$.
\end{lemma}
\proof
From Lemma \ref{Lem: interval} (a), it follows that for $r_3\in (\hat{r}_3,\bar{r}_3)$, the triangle $\triangle v_1v_2v_3$ induced by $(r_1, r_2, r_3)$ is non-degenerate. 
From (\ref{Eq: SC1}), it follows that $\varepsilon_2\kappa_1^2+\varepsilon_1\kappa_2^2 +2\kappa_1\kappa_2\eta_{12} \geq 2(\varepsilon_1\varepsilon_2+\eta_{12})\kappa_1\kappa_2 >0$. 
Combining the assumption $\Delta_{123}>0$ with (\ref{Eq: F3}) yields $\frac{\partial{h_{12,3}}}{\partial{r_3}}>0$. 
By \eqref{Eq: d}, $d_{12}$ is independent of $r_3$, and $d_{12}>0$. 
Differentiating \eqref{Eq: F6} with respect to $r_3$ gives
\begin{equation*}
\frac{\partial{\theta_{12,3}}}{\partial r_3}
=\frac{d_{12}}{d_{12}^2+(h_{12,3})^2}
\cdot\frac{\partial{h_{12,3}}}{\partial{r_3}}>0.
\end{equation*}
By Lemma \ref{Lem: theta continuous}, $\theta_{12,3}$ is a continuous function of $r_3\in [\hat{r}_3,\bar{r}_3]$. 
Therefore, $\theta_{12,3}$ is strictly increasing with respect to $r_3\in [\hat{r}_3,\bar{r}_3]$.
\qed

\begin{lemma}\label{Lem: key2}
Suppose $(S,\mathcal{T},\varepsilon,\eta)$ is a triangulated surface with a regular weight pair $(\varepsilon,\eta)$.
For any two adjacent generalized triangles $\triangle v_1v_2v_3$ and $\triangle v_1v_2v_4$ induced by a generalized label,
if $\Delta_{123}\leq 0$,
then $\theta_{12,3}+\theta_{12,4}>0$.
\end{lemma}
\proof
Since $\Delta_{123}\leq 0$, it follows from Lemma \ref{Lem: basic 1} (d) that $\eta_{12}\leq\varepsilon_1\varepsilon_2$, $\eta_{13}\leq\varepsilon_1\varepsilon_3$, $\eta_{23}\leq\varepsilon_2\varepsilon_3$, $\varepsilon_1=\varepsilon_2=\varepsilon_3=1$, and $\triangle v_1v_2v_3$ is non-degenerate. 
Therefore, $\eta_{12}\leq 1$. 

If $\eta_{12}<1$, then $h_3>0$ by \eqref{Eq: h_i}, and $h_{12,3}>0$ by (\ref{Eq: hij,k}). 
Consequently, $\theta_{12,3}>0$ by (\ref{Eq: F6}). 
Now we consider $\triangle v_1v_2v_4$. 
Assume that $\triangle v_1v_2v_4$ is non-degenerate. 
By similar arguments, we derive $h_{12,4}>0$, and then $\theta_{12,4}>0$ by \eqref{Eq: F6}. 
Assume that $\triangle v_1v_2v_4$ is degenerate. If $v_4$ is the flat vertex, 
it follows from Lemma \ref{Lem: basic 1} (c) that $\eta_{12}>\varepsilon_1\varepsilon_2=1$. 
This leads to a contradiction. 
Hence, the flat vertex of $\triangle v_1v_2v_4$ is $v_1$ or $v_2$.
From (\ref{Eq: F7}), regardless of whether the flat vertex is $v_1$ or $v_2$, we have $\theta_{12,4}=\frac{\pi}{2}$. 
Therefore, $\theta_{12,3}+\theta_{12,4}>0$. 

If $\eta_{12}=1$, then by the definition of $\Delta_{123}$ in (\ref{Eq: F1}), we derive
\begin{equation*}
0\geq \Delta_{123}
=\eta_{12}^2+\eta_{13}^2+\eta_{23}^2
+2\eta_{12}\eta_{13}\eta_{23}-1
=(\eta_{13}+\eta_{23})^2\geq 0.
\end{equation*}
This implies that $\eta_{13}=\eta_{23}=0$, and hence $\gamma_1=\gamma_2=0$ in $\triangle v_1v_2v_3$. Consequently, $h_{12,3}=0$ by (\ref{Eq: hij,k}) and $\theta_{12,3}=0$ by (\ref{Eq: F6}). 
Now we consider $\triangle v_1v_2v_4$. 
Assume that $\triangle v_1v_2v_4$ is non-degenerate. 
Then $h_{12,4}\geq0$ by (\ref{Eq: hij,k}). 
If $h_{12,4}=0$, it follows that $\gamma_1=\gamma_2=0$ in $\triangle v_1v_2v_4$, and hence $\eta_{14}=\eta_{24}=0$. 
From (\ref{Eq: SC1}), it follows that $\varepsilon_1\varepsilon_4+\eta_{14}>0$ and $\varepsilon_2\varepsilon_4+\eta_{24}>0$. 
Summing these two inequalities yields $\varepsilon_4>0$, and so $\varepsilon_4=1$. 
This contradicts the assumption that the pair $(\varepsilon,\eta)$ is regular. 
Hence, $h_{12,4}>0$ and $\theta_{12,4}>0$ by \eqref{Eq: F6}. 
Assume that $\triangle v_1v_2v_4$ is degenerate. 
By arguments analogous to those for the case $\eta_{12}<1$, the flat vertex of $\triangle v_1v_2v_4$ is $v_1$ or $v_2$,
which implies that $\theta_{12,4}=\frac{\pi}{2}$ by \eqref{Eq: F7}. 
Therefore, $\theta_{12,3}+\theta_{12,4}>0$.
\qed

\section{Maximum principle for Euclidean discrete conformal structures}
\label{Sec: MP-E}

\begin{figure}[!ht]
  \centering
  \includegraphics[scale=1]{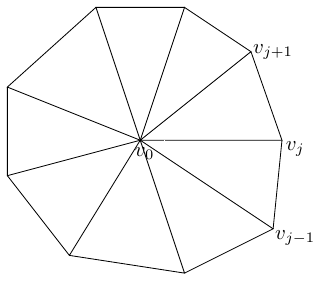}
  \caption{A star triangulation of a polygon.}
\label{figure2}
\end{figure}

In this section, we prove the maximal principle for Euclidean discrete conformal structures, i.e. Theorem \ref{Thm: MP-E}.
Recall that $P_n$ is a 1-ring neighborhood of $v_0\in V$.
See Figure \ref{figure2}. 
An assignment $r:V \to \mathbb{R}_{>0}$ can be represented as a vector in $\mathbb{R}^{n+1}$. 
Two vectors $x=(x_0,\dots, x_n)$ and $y = (y_0,\dots, y_n)$ satisfy $x\geq y$ if $x_i\geq y_i$ for all $i \in \{0,\dots, n\}$.
Theorem \ref{Thm: MP-E} can be restated as the following theorem about $r$. 

\begin{theorem}\label{Thm: MP-E2}
Let $(\varepsilon, \eta)$ be a regular weight pair on $(P_n, \mathcal{T}, \varepsilon, \eta)$.
If $r$ and $\bar{r}$ are two generalized labels on $(P_n, \mathcal{T}, \varepsilon, \eta)$ satisfying
\begin{enumerate}
\item[(a)] both $r$ and $\bar{r}$ are weighted Delaunay,
\item[(b)] the combinatorial curvatures $K_0(r)$ and $K_0(\bar{r})$ at the vertex $v_0$ satisfy $K_0(r) \leq K_0(\bar{r})$,
\item[(c)] $\max \left\{\frac{r_i}{\bar{r}_i} \mid i = 1, 2, \dots, n\right\} \leq \frac{r_0}{\bar{r}_0}$,
\end{enumerate}
then there exists a constant $c>0$ such that $r=c\bar{r}$.
\end{theorem}

To prove Theorem \ref{Thm: MP-E2}, we need the following lemma.

\begin{lemma}\label{Lem: key}
If $r, \bar{r}: \{v_0, v_1, \dots, v_n\}
\rightarrow \mathbb{R}_{>0}$ satisfies conditions (a), (b), and (c) in Theorem \ref{Thm: MP-E2}, and there exists $j \in \{1, 2, \dots, n\}$ such that $\frac{r_j}{\bar{r}_j} < \frac{r_0}{\bar{r}_0}$, then there exists $\hat{r} \in \mathbb{R}_{>0}^{n+1}$ such that
\begin{enumerate}
\item[(a)] $\hat{r}_i \geq r_i$ for $i \in \{1, \dots, n\}$,
\item[(b)] $\frac{\hat{r}_i}{\bar{r}_i} \leq \frac{\hat{r}_0}{\bar{r}_0} = \frac{r_0}{\bar{r}_0}$ for all $i = 1, 2, \dots, n$,
\item[(c)] $\hat{r}$ is a weighted Delaunay generalized label on $(P_n, \mathcal{T}, \varepsilon, \eta)$,
\item[(d)] if $\alpha(r)$ denotes the cone angle of $r$ at $v_0$, then
    \begin{equation*}
    \alpha(\hat{r}) > \alpha(r).
    \end{equation*}
\end{enumerate}
\end{lemma}
\proof
We begin by introducing some notation. 
Consider two adjacent triangles $\triangle v_0v_jv_{j\pm1}$ in $\mathcal{T}$, 
let $\theta_0^{j,j\pm 1}$ be the inner angle at the vertex $v_0$ in the triangle $\triangle v_0v_jv_{j\pm 1}$. Set
\begin{gather*}
h_j^-=h_{0j,j-1},\
h_j^+=h_{0j,j+1},\\
\theta_j^-=\theta_{0j,j-1},\
\theta_j^+=\theta_{0j,j+1},\\
\Delta_j^-=\Delta_{0j(j-1)},\
\Delta_j^+=\Delta_{0j(j+1)}.
\end{gather*}
Without loss of generality, we may assume $r_0=\bar{r}_0$ after an appropriate scaling. 
By the condition (c) in Theorem \ref{Thm: MP-E2}, it follows that $r_i\le{\bar{r}_i}$ for all $i\in\{1,2,\dots,n\}$. 
Define
\begin{gather*}
J=\{j\in\{1,2,\dots,n\}\mid r_j<{\bar{r}_j}\},\
K=\{k\in\{1,2,\dots,n\}\mid r_k={\bar{r}_k}\},\\
\gamma(r)=\sum_{j\in{J}}(\theta_j^++\theta_j^-),\
\beta(r)=\sum_{k\in{K}}(\theta_k^++\theta_k^-).
\end{gather*}

By the assumption, $J\neq\emptyset$.
By (\ref{Eq: F4}), we have
$\alpha(r)=\beta(r)+\gamma(r)$ and $\alpha(\bar{r})=\beta(\bar{r})+\gamma(\bar{r})$.
By the condition (b) in Theorem \ref{Thm: MP-E2}, it follows that $\alpha(r) \geq \alpha(\bar{r})$. 
Thus
\begin{equation}\label{Eq: F5}
\beta(r)+\gamma(r)\geq\beta(\bar{r})+\gamma(\bar{r}).
\end{equation}

\textbf{Claim 1}: For any $j\in J$, we have $\theta_{0}^{j-1,j}(r)<\pi$ and $\theta_{0}^{j,j+1}(r)<\pi$.

It suffices to show that for any $j\in J$, if the triangle $\triangle v_0v_{j}v_{j-1}$ is degenerate, then $v_0$ cannot be the flat vertex. 
If not, assume that there exists some $j\in J$ such that $v_0$ is the flat vertex of the degenerate triangle $\triangle v_0v_{j}v_{j-1}$ induced by $(r_0,r_j,r_{j-1})$ via (\ref{Eq: length-E2}). 
By Lemma \ref{Lem: basic 2} (c), $\kappa_0$ is a function $\kappa_0=f(\kappa_{j-1},\kappa_{j})$ of $\kappa_{j-1}$ and $\kappa_{j}$ in this case. Furthermore, the function $\kappa_0=f(\kappa_{j-1},\kappa_{j})$ is strictly increasing in both $\kappa_j$ and $\kappa_{j-1}$. 
Since $\kappa_j>\bar{\kappa}_j$ and $\kappa_{j-1}\geq\bar{\kappa}_{j-1}$,
it follows that
\begin{equation*}
\bar{\kappa}_0=\kappa_0=f(\kappa_{j-1},\kappa_{j})
>f(\bar{\kappa}_{j-1},\bar{\kappa}_{j}).
\end{equation*}
By Lemma \ref{Lem: basic 2} (a), $(\bar{r}_0, \bar{r}_j, \bar{r}_{j-1})\in V_0\setminus \partial V_0$.
Applying Lemma \ref{Lem: basic 2} (c) again, we conclude that $(\bar{r}_0,\bar{r}_j,\bar{r}_{j-1})$ is not a generalized label. 
This contradicts the assumption in Theorem \ref{Thm: MP-E2} that $\bar{r}$ is a generalized label.

\textbf{Claim 2}: There exists $j\in J$ such that $\theta_j^+(r)+\theta_j^-(r)>0$.

If there exists some $j\in J$ such that $\Delta_j^+\leq0$ or $\Delta_j^-\leq0$,
it follows that $\theta_j^+(r)+\theta_j^-(r)>0$ by Lemma \ref{Lem: key2}. 
Therefore, we may assume that for all $j\in J$, $\Delta_j^+>0$ and $\Delta_j^->0$.

If $K\ne\emptyset$, then there exists $i\in K$ such that either $i-1\in J$ or $i+1\in J$. 
Without loss of generality, we may assume $i-1\in J$. 
Consequently, $\Delta_{i-1}^+=\Delta_i^->0$. 
By Lemma \ref{Lem: monotonicity}, the function $\theta_i^-$ is strictly increasing in $r_{i-1}$. 
It follows that $\theta_i^-(r)<\theta_i^-(\bar{r})$. 
If $i+1\in K$, then $\theta_i^+(r)=\theta_i^+(\bar{r})$. 
If $i+1\in J$, then $\Delta_{i-1}^-=\Delta_i^+>0$. 
By Lemma \ref{Lem: monotonicity}, the function $\theta_i^+$ is strictly increasing in $r_{i+1}$. 
It follows that $\theta_i^+(r)<\theta_i^+(\bar{r})$. Therefore, $\theta_i^-(r)+\theta_i^+(r)<\theta_i^-(\bar{r})+\theta_i^+(\bar{r})$.
By induction, it follows that $\beta(r)<\beta(\bar{r})$. 
Consequently, by (\ref{Eq: F5}), we have $0\leq\gamma(\bar{r})<\gamma(r)$. 
By the definition of $\gamma(r)$, there exists $j\in J$ such that $\theta_j^+(r)+\theta_j^-(r)>0$.

If $K=\emptyset$, then $J=\{1,\dots,n\}$, and
\begin{equation*}
\gamma(r)=\sum_{j\in{J}}(\theta_j^+(r)+\theta_j^-(r))
=\alpha(r)\ge0.
\end{equation*}
If $\alpha(r)>0$, then there exists $j\in J$ such that $\theta_j^+(r)+\theta_j^-(r)>0$. 
We will show that the case $\alpha(r)=0$ is impossible.

If $\alpha(r)=0$, then for every triangle $\triangle v_0v_jv_{j-1}$, it holds that $\theta_0^{j,j-1}=0$ for all $j\in \{1,\dots,n\}$. 
Consequently, all triangles are degenerate. 
For every triangle $\triangle v_0v_jv_{j-1}$, the flat vertex is either $v_j$ or $v_{j-1}$. 
By (\ref{Eq: F7}), it follows that $\{\theta_j^-(r),\theta_{j-1}^+(r)\}=\{\frac{\pi}{2}, -\frac{\pi}{2}\}$ for all $j\in \{1, \dots, n\}$. 
Without loss of generality, we may assume that $v_1$ is the flat vertex of $\triangle v_0v_1v_2$. 
It follows from (\ref{Eq: F7}) that $\theta_1^+(r)=\frac{\pi}{2}$ and $\theta_2^-(r)=-\frac{\pi}{2}$.
Moreover, $l_{02}(r)=l_{01}(r)+l_{12}(r)>l_{01}(r)$.
By the weighted Delaunay condition in (\ref{Eq: F8}), it follows that $\theta_2^+(r)=\frac{\pi}{2}$.
Consequently, $\theta_3^-(r)=-\frac{\pi}{2}$ due to $\theta_0^{23}=0$, and $v_2$ is the flat vertex of $\triangle v_0v_2v_3$ by (\ref{Eq: F7}). 
Therefore, $l_{03}(r)=l_{02}(r)+l_{23}(r)>l_{02}(r)$. By induction, we obtain a contradiction
\begin{equation*}
l_{01}(r)<l_{02}(r)<\dots<l_{0n}(r)<l_{01}(r).
\end{equation*}
This concludes the proof of Claim 2.


For a fixed index $j \in J$ as stated in Claim 2, we have  $\theta_j^+(r)+\theta_j^-(r)>0$.

\textbf{Claim 3}:
Assume that $\triangle v_0v_jv_{j-1}$ is degenerate.
Then the flat vertex must be $v_j$.

According to Claim 1, $v_0$ cannot be the flat vertex.
If $v_{j-1}$ is the flat vertex, then by (\ref{Eq: F7}), we obtain $\theta_j^-=-\frac{\pi}{2}$.
Since $\theta_j^+ \in \left[-\frac{\pi}{2}, \frac{\pi}{2}\right]$,
it follows that $\theta_j^+(r) + \theta_j^-(r) \leq 0$.
This leads to a contradiction.

In the following, we prove that there exists $\epsilon>0$ such that the vector $\hat{r}=(r_0,\dots,r_j+t,\dots,r_n)$ satisfies Lemma \ref{Lem: key} for all $t\in(0,\epsilon)$. 
It is straightforward to verify that for $t\in(0,\bar{r}_j-r_j)$, $\hat{r}$  satisfies the conclusions (a) and (b) of Lemma \ref{Lem: key}.

To see the conclusion (c) of Lemma \ref{Lem: key}, we first prove that there exists $\epsilon>0$ such that
$\hat{r}$ is a generalized label on $(P_n, \mathcal{T}, \varepsilon, \eta)$ for $t\in (0,\epsilon)$.
Then we show that $\hat{r}$ satisfies the weighted Delaunay condition.

If $\triangle v_0v_jv_{j-1}$ induced by $r$ is non-degenerate,  
then $\triangle v_0v_jv_{j-1}$ induced by $\hat{r}$ remains non-degenerate due to continuity.  
If $\triangle v_0v_jv_{j-1}$ induced by $r$ is degenerate, then by Claim 3, the flat vertex is $v_j$.  
By Lemma \ref{Lem: interval} (b), there exists $\epsilon>0$ such that for all $t\in(0,\epsilon)$,  
$\triangle v_0v_jv_{j-1}$ induced by $\hat{r}$ is non-degenerate.  
A similar argument applies to the triangle $\triangle v_0v_jv_{j+1}$.  
Therefore, there exists $\epsilon>0$ such that for all $t\in(0,\epsilon)$,  
$\hat{r}$ is a generalized label on $(P_n,\mathcal{T},\varepsilon,\eta)$.  
Moreover, the triangles $\triangle v_0v_jv_{j\pm1}$ induced by $\hat{r}$ are non-degenerate.

Since $\hat{r}$ differs from $r$ only at the $j$-th position,  
it suffices to analyze the edges $v_0v_j$ and $v_0v_{j\pm1}$.  
For the edge $v_0v_j$, by combining $\theta_j^+(r)+\theta_j^-(r)>0$ with the continuity of $\theta_j^{\pm}$ in Lemma \ref{Lem: theta continuous},
it follows that $\theta_j^+(\hat{r})+\theta_j^-(\hat{r})>0$ for sufficiently small $t>0$. 
For the edge $v_0v_{j-1}$, it holds that $\theta_{j-1}^-(r)=\theta_{j-1}^-(\hat{r})$.  
If $\Delta_j^->0$, then it follows from Lemma \ref{Lem: monotonicity} that $\theta_{j-1}^+(r)<\theta_{j-1}^+(\hat{r})$ for all $t\in(0,\bar{r}_j-r_j)$.  
This implies that $\theta_{j-1}^+(\hat{r})+\theta_{j-1}^-(\hat{r})  >\theta_{j-1}^+(r)+\theta_{j-1}^-(r)\ge0$.  
If $\Delta_j^-\leq0$, then it follows from Lemma \ref{Lem: key2} that $\theta_{j-1}^+(r)+\theta_{j-1}^-(r)>0$.  
The conclusion is derived from the continuity of $\theta_{j-1}^{\pm}$ in Lemma \ref{Lem: theta continuous}.  
Therefore, there exists $\epsilon>0$ such that the edge $v_0v_{j-1}$ satisfies the weighted Delaunay condition in $\hat{r}$ for all $t\in(0,\epsilon)$.  
A similar argument applies to the edge $v_0v_{j+1}$.

To see the conclusion (d) of Lemma \ref{Lem: key}, following the arguments presented in part (c), there exists $\epsilon>0$ such that the triangles $\triangle v_0v_jv_{j\pm1}$ are non-degenerate in $\hat{r}$ and $\theta_{j}^+(\hat{r})+\theta_{j}^-(\hat{r})>0$ for all $t\in(0,\epsilon)$.  
By Corollary \ref{Cor: weighted Delaunay 2},  
it follows that $h_j^+(\hat{r})+h_j^-(\hat{r})>0$ for all $t\in(0,\epsilon)$.  
Note that the function $\alpha(\hat{r})$ is continuous for $t\in[0,\epsilon]$ and smooth for $t\in(0,\epsilon)$.  
By Lemma \ref{Lem: basic 1} (b), we derive 
\begin{equation*}  
\frac{\partial \alpha}{\partial t}(\hat{r})  
=\frac{h_j^+(\hat{r})+h_j^-(\hat{r})}{l_{0j}}>0,\quad  
t\in(0,\epsilon).  
\end{equation*}  
Consequently, $\alpha(\hat{r})>\alpha(r)$ for all $t\in(0,\epsilon)$.
\qed

Now we can prove Theorem \ref{Thm: MP-E2}.

\noindent\textbf{Proof for Theorem \ref{Thm: MP-E2}:}
The proof is similar to that of Theorem 3.1 in \cite{LSW}.  
Without loss of generality, we perform a scaling transformation to assume that $r_0=\bar{r}_0$ and $r_i\le\bar{r}_i$ for all $i=1,2,\dots,n$.  
We proceed by contradiction. Suppose Theorem \ref{Thm: MP-E2} does not hold.  
Then there exists a vector $r$ satisfying $r_0=\bar{r}_0$ and $r_i\le\bar{r}_i$ for all $i=1,2,\dots,n$, with at least one index $i_0$ such that $r_{i_0}<\bar{r}_{i_0}$ and $\alpha(\bar{r})\leq\alpha(r)$.  
By Lemma \ref{Lem: key}, after replacing $r$ with $\hat{r}$, we can further deduce that
\begin{equation*}  
\alpha(\bar{r})<\alpha(r).  
\end{equation*}  
Define the set  
\begin{equation*}  
X:=\{x\in\mathbb{R}^{n+1}\mid r\le x\le \bar{r}, \text{$x$ is a weighted Delaunay generalized label on $(P_n, \mathcal{T}, \varepsilon, \eta)$}\}.  
\end{equation*}  
Clearly, $r\in X$, and $X$ is bounded.  
By Lemma \ref{Lem: theta continuous}, the set $X$ is closed in $\mathbb{R}^{n+1}$.  
Therefore, $X$ is compact, and the function $\alpha(x)$ attains its maximum on $X$. 
Let $t\in X$ be a point where the continuous function $\alpha(x)$ achieves its maximum on $X$.  
If $t\neq\bar{r}$, then by Lemma \ref{Lem: key},  
there exists a weighted Delaunay generalized label $\hat{t}$ on $(P_n, \mathcal{T}, \varepsilon, \eta)$ such that $\hat{t}\ge t$, $\hat{t}_0=\bar{r}_0$, $\hat{t}\le\bar{r}$, and $\alpha(\hat{t})>\alpha(t)$.
This contradicts the assumption that $t$ is the maximum point of $\alpha(x)$.  
Therefore, $t=\bar{r}$, and 
\begin{equation*}  
\alpha(\bar{r})=\alpha(t)\ge\alpha(r)>\alpha(\bar{r}).  
\end{equation*}  
This leads to a contradiction. 
\qed

\begin{remark}
The maximum principle in Theorem \ref{Thm: MP-E2} unifies and generalizes a wide range of discrete maximum principles on surfaces.  
Specifically, when $\varepsilon\equiv0$ and $\eta: E\rightarrow \mathbb{R}_{>0}$,  
Theorem \ref{Thm: MP-E2} reduces to the maximum principle for Luo's vertex scalings established by Luo-Sun-Wu \cite{LSW} and Dai-Ge-Ma \cite{DGM}.  
When $\varepsilon\equiv1$ and $\eta: E\rightarrow \mathbb{R}_{\geq0}$,  
Theorem \ref{Thm: MP-E2} reduces to the maximum principle for inversive distance circle packing established by Luo-Xu-Zhang \cite{LXZ}.  
Note that the maximal principle for the case $\varepsilon\equiv1$ and $\eta: E\rightarrow (-1,1]$ has been established in \cite{LXZ}.
\end{remark}

\begin{remark}
In the case where $\varepsilon\equiv1$, Luo-Xu-Zhang \cite{LXZ} provided a counterexample demonstrating that the discrete maximum principle fails when the weight $\eta$ simultaneously takes values in both $(-1,0)$ and $(1,+\infty)$ on a triangulated surface.  
For further details, see Remark 3.3 in \cite{LXZ}. 
We now investigate the remaining case where $\varepsilon: V\rightarrow \{0,1\}$ (with $\varepsilon\not\equiv0$ and $\varepsilon\not\equiv1$) and $\eta: E\rightarrow (-1,1]$ satisfy the structure conditions (\ref{Eq: SC1}) and (\ref{Eq: SC2}).  
In the following, we present a counterexample to demonstrate that Theorem \ref{Thm: MP-E2} does not hold in this case.  
Consider a polygon disk $P_4$ with four boundary vertices $v_1, v_2, v_3, v_4$ and a unique interior vertex $v_0$.  
Define $\varepsilon_0=\varepsilon_1=\varepsilon_3=1$,  
$\varepsilon_2=\varepsilon_4=0$,  
$\eta_{01}=\eta_{03}=-\frac{1}{2}$,  
$\eta_{02}=\eta_{04}=\frac{1}{2}$,  
and $\eta_{12}=\eta_{23}=\eta_{34}=\eta_{14}=\frac{1}{4}$.  
It is straightforward to check that the pair $(\eta,\varepsilon)$ is regular and satisfies the structure conditions (\ref{Eq: SC1}) and (\ref{Eq: SC2}).
We further define $r_0=1$, $r_1=r_3=2$, and $r_2=r_4=c^{-1}$ for $c>0$.  
For the triangle $\triangle v_0v_1v_2$, we have $Q=\frac{3}{4}c(c+1)>0$,  
which implies that $\triangle v_0v_1v_2$ is non-degenerate.  
Direct calculations yield  
$\gamma_0=\varepsilon_0\eta_{12}+\eta_{01}\eta_{02}=0$,  
$\gamma_1=\varepsilon_1\eta_{02}+\eta_{01}\eta_{12}=\frac{3}{8}$,  
and $\gamma_2=\varepsilon_2\eta_{01}+\eta_{02}\eta_{12}=\frac{1}{8}$.  
Furthermore, $h_{01,2}>0$ and $h_{02,1}=0$.  
Moreover, $l_{01}^2+l_{02}^2=l_{12}^2$,  
which implies that $\triangle v_0v_1v_2$ is a right triangle with $\angle v_1v_0v_2=\frac{\pi}{2}$.  
It is straightforward to verify that all triangles are congruent.  
Therefore,  
$h_{01,2}=h_{01,4}=h_{03,2}=h_{03,4}>0$ and  
$h_{02,1}=h_{02,3}=h_{04,1}=h_{04,3}=0$.  
Thus,  
$r=(1,2,c^{-1},2,c^{-1})\in \mathbb{R}^5_{>0}$ is a weighted Delaunay admissible label on $P_4$ for any $c>0$.  
The cone angle $\alpha(r)$ at the vertex $v_0$ is $2\pi$ for any $c>0$.  
This implies that the discrete maximum principle does not hold in this case.
\end{remark}

\section{Maximum principle for hyperbolic discrete conformal structures}
\label{Sec: MP-H}

\begin{figure}[!ht]
\centering
\includegraphics[scale=0.7]{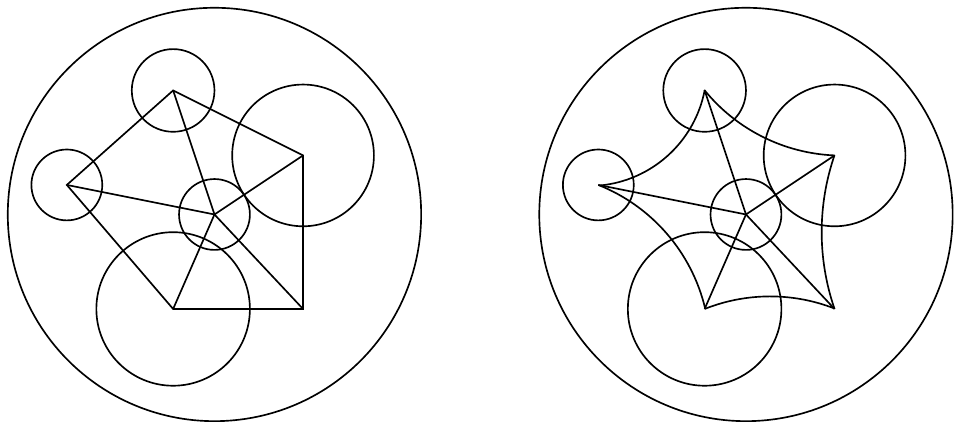}
\caption{The PL and PH metric from one configuration of vertex-circles.}
\label{figure3}
\end{figure}

In this section, we prove the maximal principle for hyperbolic discrete conformal structures from its Euclidean counterparts. We will establish the relationships between the Euclidean discrete conformal factors and the hyperbolic discrete conformal factors, using the configuration in Figure \ref{figure3} as a model. 

Recall that $P_n$ is an $n$-sided star triangulation of a polygon. In the configuration in Figure \ref{figure3}, a vertex-circle is located at the origin $v_0$, together with a finite collection of vertex-circles corresponding to $v_1, \cdots, v_n$ of $P_n$ inside the unit disk $\mathbb{D}$. Then there is a PL metric induced by a label on $(P_n, \mathcal{T}, \varepsilon, \eta^e)$ by connecting the Euclidean centers of the vertex-circles by Euclidean segments. 
Meanwhile, there is also a PH metric induced by another label on $(P_n, \mathcal{T}, \varepsilon, \eta^h)$ by connecting the hyperbolic centers of the vertex-circles by hyperbolic segments in $\mathbb{D}$. 
It is found that, if the Euclidean and hyperbolic labels are chosen appropriately, then they could induce the same vertex circles in Figure \ref{figure3} with $\eta^e=\eta^h$.
Note that the locations of the Euclidean center and the hyperbolic center of a vertex-circle are different in general, except the vertex-circle at the vertex $v_0$ located at the origin. 
Moreover,  the induced PL metric is weighted Delaunay if and only if the induced PH metric is weighted Delaunay, 
since the intersections angles of circles are conformally invariant. 
Note that for the configuration in Figure \ref{figure3}, the inner angles at the origin $v_0$ in the Euclidean background geometry are the same
as the corresponding inner angles at $v_0$ in the hyperbolic background geometry. As a result, the combinatorial curvatures of the PL metric and the PH metric at $v_0$ are identical. 

To distinguish between the Euclidean and hyperbolic cases and to handle the general cases, we adopt the convention of using capital letters for Euclidean quantities and lowercase letters for hyperbolic quantities.
\begin{enumerate}
    \item In the Euclidean case, the labels are denoted as $U_i=F_i$ and $e^{F_i}=R_i$, where $\varepsilon$ takes values in $\{0,1\}$.
    \item In the hyperbolic case, the labels are given by $e^{f_i}=\sinh r_i$, and $e^{u_i}=\tanh\frac{r_i}{2}$ when $\varepsilon_i=1$, and $u_j=f_j$ when $\varepsilon_j=0$.
\end{enumerate}
The parameter $\varepsilon$ remains identical at the same vertex in PL metrics and PH metrics.
The equations (\ref{Eq: length-E}) and (\ref{Eq: length-E2}) for the Euclidean discrete conformal structures can be reformulated as 
\begin{equation}\label{Eq: length-E3}
\begin{aligned}
L^2_{ij}
&=\varepsilon_ie^{2F_i}+\varepsilon_je^{2F_j}
+2\eta^e_{ij}e^{F_i+F_j}\\
&=\varepsilon_ie^{2U_i}+\varepsilon_je^{2U_j}
+2\eta^e_{ij}e^{U_i+U_j}\\
&=\varepsilon_iR_i^2+\varepsilon_jR_j^2+2\eta^e_{ij}R_iR_j.
\end{aligned}
\end{equation}

In the rest of this section, we will first give the relationships between hyperbolic and Euclidean labels. 
We will also show that the hyperbolic weights are identical with the Euclidean weights under these relationships. 
To prove a hyperbolic maximal principle, we apply a similar transformation to the configuration in the plane. This transformation induces variations in both the Euclidean and hyperbolic labels.  
Finally, we examine the relationships between the hyperbolic and Euclidean discrete conformal factors.

\subsection{Relationships between hyperbolic and Euclidean labels}\label{Sec: key}

\begin{figure}[!ht]
\centering
\includegraphics[scale=0.7]{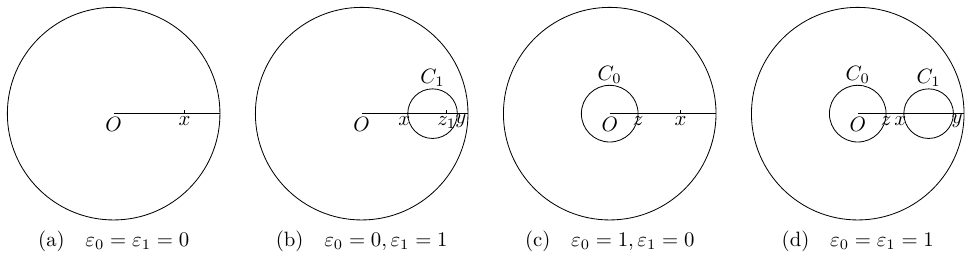}
\caption{Four cases of the configurations. Note that the exterior intersection angle of circles $C_0$ and $C_1$ is in $[0,\frac{\pi}{2}]$.}
\label{figure1}
\end{figure}

Let $e=v_0v_1$ be a geodesic segment in $\mathbb{D}$, where the endpoint $v_0$ coincides with the origin $O$ of $\mathbb{D}$, and the endpoint $v_1$ lies on the positive $x$-axis.
We classify the configurations into the following four distinct cases:
\begin{enumerate}
   \item [(a)] $\varepsilon_0=\varepsilon_1=0$;
 \item  [(b)] $\varepsilon_0=0,\varepsilon_1=1$;
\item  [(c)] $\varepsilon_0=1,\varepsilon_1=0$;
\item [(d)] $\varepsilon_0=\varepsilon_1=1$.
\end{enumerate}
As discussed in Subsections \ref{Sec: 1.1} and \ref{Sec: 1.2}, 
if $\varepsilon_i=1$, then there exists a vertex-circle $C_i$ with hyperbolic radius $r_i$ and Euclidean radius $R_i$;
if $\varepsilon_i=0$, then we regard $C_i$ as shrinking into a single point.
See Figure \ref{figure1}.

In the Poincar\'{e} disk model $\mathbb{D}$, the hyperbolic distance $d(z_1,z_2)$ between any two points $z_1$ and $z_2$ is given by
\begin{equation}\label{Eq: HD1}
\cosh d(z_1,z_2)
=1+\frac{2|z_1-z_2|^2}{(1-|z_1|^2)(1-|z_2|^2)},
\end{equation}
or equivalently,
\begin{equation}\label{Eq: HD2}
\sinh \frac{d(z_1,z_2)}{2}
=\frac{|z_1-z_2|}{\sqrt{(1-|z_1|^2)(1-|z_2|^2)}}.
\end{equation}
For a detailed derivation of (\ref{Eq: HD1}), please refer to \cite[Theorem 4.5.1]{Ratcliffe}.

\textbf{Case 1}: 
Assume that $\varepsilon_0=\varepsilon_1=0$.
See Figure \ref{figure1} (a).
Let $x$ denote a real number representing the position of $v_1$, where $0<x<1$.
From (\ref{Eq: length-E3}), it follows that
\begin{equation}\label{Eq: case 1 a} 
x^2=L_{01}^2=2\eta_{01}^ee^{U_0+U_1}.
\end{equation}
By combining (\ref{Eq: HD1}) and (\ref{Eq: length-H}), we obtain
\begin{equation*}
1+\frac{2x^2}{1-x^2}=\cosh l_{01}
=1+\eta_{01}^he^{u_0+u_1},
\end{equation*}
which is equivalent to
\begin{equation}\label{Eq: case 1 b} 
\frac{2x^2}{1-x^2}
=\eta_{01}^he^{u_0+u_1}.
\end{equation}
Combining (\ref{Eq: case 1 a}) and (\ref{Eq: case 1 b}), we derive
\begin{equation}\label{Eq: F17}
\frac{\eta_{01}^h}{\eta_{01}^e}\cdot e^{u_0-U_0}\cdot e^{u_1-U_1}=\frac{4}{1-x^2}.
\end{equation}
To ensure $\eta_{01}^e=\eta_{01}^h$, we can set
\begin{equation*}
u_0=U_0+\ln 2
\end{equation*}
and
\begin{equation*}
u_1=U_1+\ln\frac{2}{1-x^2}.
\end{equation*}

\textbf{Case 2}:
Assume that $\varepsilon_0=0$ and $\varepsilon_1=1$.
See Figure \ref{figure1} (b).
There exists a circle $C_1$ with its hyperbolic center at $v_1$.
The Euclidean and the hyperbolic radii of $C_1$ are denoted by $R_1$ and $r_1$, respectively.
Let $z_1$ be the real number representing the position of the hyperbolic center $v_1$.
The intersection points of $C_1$ with the $x$-axis are denoted by $x$ and $y$, and thus $0<x<z_1<y<1$.

From the relation $e^{U_1}=R_1=\frac{y-x}{2}$, it is straightforward to deduce that
\begin{equation}\label{Eq: case 2 a}
U_1=\ln \frac{y-x}{2}.
\end{equation}
Substituting \eqref{Eq: case 2 a} into (\ref{Eq: length-E3}) yields
\begin{equation*}
\left(\frac{x+y}{2}\right)^2=L_{01}^2
=\left(\frac{y-x}{2}\right)^2+2\eta_{01}^ee^{U_0}\cdot \frac{y-x}{2}.
\end{equation*}
This implies that
\begin{equation}\label{Eq: Euclidean i0j1}
\eta_{01}^ee^{U_0}=\frac{xy}{y-x}.
\end{equation}
From (\ref{Eq: HD2}), it follows that
\begin{equation}\label{Eq: F45}
\sinh r_1=\frac{y-x}{\sqrt{1-x^2}\sqrt{1-y^2}}
\quad \text{and} \quad
\cosh r_1=\frac{1-xy}{\sqrt{1-x^2}\sqrt{1-y^2}}.
\end{equation}
Consequently,
\begin{equation}\label{Eq: F49}
e^{u_1}=\tanh \frac{r_1}{2}
=\frac{\sinh r_1}{1+\cosh r_1}
=\frac{y-x}{1-xy+\sqrt{(1-x^2)(1-y^2)}}.
\end{equation}
Therefore,
\begin{equation}\label{Eq: F19}
\begin{aligned}
u_1=&\ln \frac{y-x}{2}+\ln \frac{2}{1-xy+\sqrt{(1-x^2)(1-y^2)}}\\
=&U_1+\ln \frac{2}{1-xy+\sqrt{(1-x^2)(1-y^2)}},
\end{aligned}
\end{equation}
where the equation (\ref{Eq: case 2 a}) is used in the last line.

By combining (\ref{Eq: length-H}) and (\ref{Eq: F45}), we obtain
\begin{equation}\label{Eq: F46}
\begin{aligned}
\cosh l_{01}
&=\sqrt{1+e^{2f_1}}
+\eta_{01}^he^{f_0+f_1}\\
&=\cosh r_1+\eta_{01}^he^{u_0}\sinh r_1\\
&=\frac{(1-xy)+\eta_{01}^he^{u_0}(y-x)}{\sqrt{1-x^2}\sqrt{1-y^2}}.
\end{aligned}
\end{equation}
Furthermore, by (\ref{Eq: HD2}), we have
\begin{equation*}
\frac{z_1-x}{\sqrt{1-x^2}\sqrt{1-z_1^2}}
=\sinh\frac{d(x,z_1)}{2}
=\sinh\frac{d(y,z_1)}{2}
=\frac{y-z_1}{\sqrt{1-y^2}\sqrt{1-z_1^2}}.
\end{equation*}
This implies that
\begin{equation}\label{Eq: F42}
z_1
=\frac{1+xy-\sqrt{1-x^2}\sqrt{1-y^2}}{x+y}
=\frac{x+y}{1+xy+\sqrt{1-x^2}\sqrt{1-y^2}}.
\end{equation}
Consequently,
\begin{equation}\label{Eq: F43}
1+z^2_1=\frac{2(1+xy)}{1+xy+\sqrt{1-x^2}\sqrt{1-y^2}}
\quad \text{and} \quad
1-z^2_1=\frac{2\sqrt{1-x^2}\sqrt{1-y^2}}{1+xy+\sqrt{1-x^2}\sqrt{1-y^2}}.
\end{equation}
Combining (\ref{Eq: F43}) with (\ref{Eq: HD1}) yields
\begin{equation}\label{Eq: F37}
\cosh l_{01}=\frac{1+z_1^2}{1-z_1^2}
=\frac{1+xy}{\sqrt{1-x^2}\sqrt{1-y^2}}.
\end{equation}
By comparing (\ref{Eq: F46}) and (\ref{Eq: F37}), we derive
\begin{equation}\label{Eq: hyperbolic i0j1}
\eta_{01}^he^{u_0}=\frac{2xy}{y-x}.
\end{equation}
By combining (\ref{Eq: Euclidean i0j1}) and (\ref{Eq: hyperbolic i0j1}),
to ensure $\eta_{01}^e=\eta_{01}^h$, we have
\begin{equation*}
u_0=U_0+\ln 2.
\end{equation*}
This justifies the choice of $u_0=U_0+\ln 2$ in Case 1.

\textbf{Case 3}:
Assume that $\varepsilon_0=1$ and $\varepsilon_1=0$.
See Figure \ref{figure1} (c).
There exists a circle $C_0$ centered at the origin.
The Euclidean and hyperbolic radii of $C_0$ are denoted by $R_0$ and $r_0$, respectively.
The intersection point of $C_0$ with the positive $x$-axis is denoted by $z$, and hence $R_0=z$.
Let $x$ denote a real number representing the position of $v_1$, and then $0<z<x<1$.

Since $e^{U_0}=R_0=z$, it follows from (\ref{Eq: length-E3}) that
\begin{equation*}
x^2=L_{01}^2
=e^{2U_0}+2\eta_{01}^ee^{U_0+U_1}
=z^2+2z\cdot \eta_{01}^ee^{U_1}.
\end{equation*}
This implies that
\begin{equation}\label{Eq: Euclidean i1j0}
\eta_{01}^ee^{U_1}=\frac{x^2-z^2}{2z}.
\end{equation}
By (\ref{Eq: HD1}), we have
\begin{equation}\label{Eq: F47}
\cosh r_0=\frac{1+z^2}{1-z^2}
\quad \text{and} \quad
\sinh r_0=\frac{2z}{1-z^2}.
\end{equation}
Consequently,
\begin{equation*}
e^{u_0}=\tanh \frac{r_0}{2}
=\frac{\sinh r_0}{1+\cosh r_0}
=z.
\end{equation*}
Therefore,
\begin{equation*}
u_0=U_0.
\end{equation*}
By combining (\ref{Eq: HD1}) and (\ref{Eq: length-H}), we obtain
\begin{equation}\label{Eq: F48}
\frac{1+x^2}{1-x^2}=\cosh l_{01}
=\sqrt{1+e^{2f_0}}
+\eta_{01}^he^{f_0+f_1}
=\cosh r_0+\eta_{01}^he^{u_1}\sinh r_0.
\end{equation}
Substituting (\ref{Eq: F47}) into (\ref{Eq: F48}) yields
\begin{equation}\label{Eq: hyperbolic i1j0}
\eta_{01}^he^{u_1}=\frac{x^2-z^2}{z(1-x^2)}.
\end{equation}
By combining (\ref{Eq: Euclidean i1j0}) and (\ref{Eq: hyperbolic i1j0}),
to ensure $\eta_{01}^e=\eta_{01}^h$, we obtain
\begin{equation*}
u_1=U_1+\ln \frac{2}{1-x^2}.
\end{equation*}

\textbf{Case 4}:
Assume that $\varepsilon_0=\varepsilon_1=1$.
See Figure \ref{figure1} (d).
The notations for circles $C_0$ and $C_1$ are consistent with those in Case 2 and Case 3, respectively.
From Case 3, it follows that
\begin{equation*}
u_0=U_0.
\end{equation*}
From Case 2, it follows that
\begin{equation*}
u_1=U_1+\ln \frac{2}{1-xy+\sqrt{(1-x^2)(1-y^2)}}.
\end{equation*}
It is straightforward to check that $\eta_{01}^e=\eta_{01}^h$ in this case.


Based on the preceding analysis, we establish the following conventions.

\begin{convention}\label{Convention}
Let $e=v_iv_j$ denote a geodesic segment in $\mathbb{D}$ with $v_i=O$.
\begin{enumerate}
\item [(a)]
If $\varepsilon_i=0$, then we set $u_i=U_i+\ln 2$;
\item [(b)]
If $\varepsilon_i=1$, then we set $u_i=U_i$;
\item [(c)]
If $\varepsilon_j=0$, then we set $u_j=U_j+\ln\frac{2}{1-|z_j|^2}$,
where $z_j$ is a point representing $v_j$;
\item [(d)]
If $\varepsilon_j=1$, then we set $u_j=U_j+\ln \frac{2}{1-xy+\sqrt{(1-x^2)(1-y^2)}}$,
where $x$ and $y$ denote the real numbers representing the minimal and maximal Euclidean distances, respectively, from the points of $C_j$ (with hyperbolic center $v_j$) to the origin $O$.
\end{enumerate}
\end{convention}

\begin{figure}[!ht]
\centering
\includegraphics[scale=0.9]{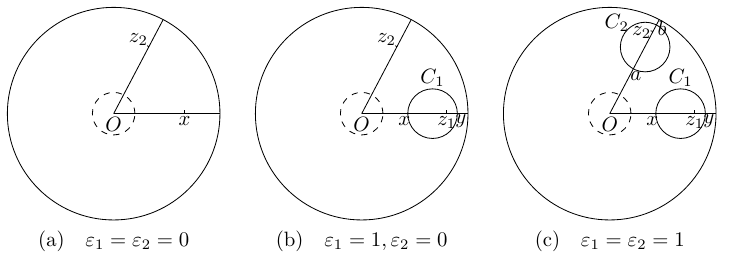}
\caption{Three cases. Note that the dashed circle centered at $O$ represents $\varepsilon_0\in \{0,1\}$.}
\label{figure4}
\end{figure}

In the subsequent analysis, we will prove that $\eta^h=\eta^e$ for a triangle $\triangle v_0v_1v_2$ under Convention \ref{Convention}.
The radial direction has already been verified in the four cases discussed above.
It remains to prove that $\eta^h_{12}=\eta^e_{12}$.
By applying a rotation transformation, we may assume that the vertex $v_1$ lies on the positive $x$-axis.
We now consider the following three distinct cases: 
\begin{enumerate}
    \item [(I)] $\varepsilon_1=\varepsilon_2=0$; 
    \item [(II)] $\varepsilon_1=1$, $\varepsilon_2=0$; 
    \item [(III)] $\varepsilon_1=1$, $\varepsilon_2=1$.
\end{enumerate}
See Figure \ref{figure4}.

\textbf{Case I}:
Assume that $\varepsilon_1=\varepsilon_2=0$.
See Figure \ref{figure4} (a).
Let $x$ denote a real number representing the position of $v_1$, and let $z_2$ denote a complex number representing the position of $v_2$.
According to Convention \ref{Convention}, we obtain
\begin{equation}\label{case1a}
u_1=U_1+\ln\frac{2}{1-x^2}
 \quad \text{and} \quad
u_2=U_2+\ln\frac{2}{1-|z_2|^2}.
\end{equation}
From (\ref{Eq: length-E3}), it follows that
\begin{equation}\label{case1b}
|z_2-x|^2=L_{12}^2
=2\eta_{12}^ee^{U_1+U_2}.
\end{equation}
Combining (\ref{Eq: HD1}) and (\ref{Eq: length-H}) yields
\begin{equation}\label{case1c}
1+\frac{2|z_2-x|^2}{(1-x^2)(1-|z_2|^2)}=\cosh l_{12}
=1+\eta_{12}^he^{u_1+u_2}.
\end{equation}
By combining (\ref{case1a}), (\ref{case1b}), and (\ref{case1c}), we conclude that $\eta_{12}^h=\eta_{12}^e$.

\textbf{Case II}:
Assume that $\varepsilon_1=1$ and $\varepsilon_2=0$.
See Figure \ref{figure4} (b).
The notations for the circle $C_1$ are consistent with those in Case 2.
Let $z_2$ denote a complex number representing the position of $v_2$.
Since $e^{U_1}=R_1=\frac{y-x}{2}$, it follows from (\ref{Eq: length-E3}) that
\begin{equation*}
\left|z_2-\frac{x+y}{2}\right|^2=L_{12}^2
=\left(\frac{y-x}{2}\right)^2+2\eta_{12}^ee^{U_2}\cdot \frac{y-x}{2},
\end{equation*}
which implies
\begin{equation}\label{Eq: F50}
|z_2|^2-\frac{x+y}{2}(z_2+\bar{z}_2)+xy
=(y-x)\eta_{12}^ee^{U_2}.
\end{equation}
Combining (\ref{Eq: length-H}) and (\ref{Eq: F45}) gives
\begin{equation}\label{CaseIIa}
\begin{aligned}
\cosh l_{12}
&=\sqrt{1+e^{2f_1}}+\eta_{12}^he^{f_1+f_2}\\
&=\cosh r_1+\eta_{12}^he^{u_2}\sinh r_1\\
&=\frac{(1-xy)+\eta^h_{12}e^{u_2}(y-x)}{\sqrt{1-x^2}\sqrt{1-y^2}},
\end{aligned}
\end{equation}
Furthermore, by (\ref{Eq: HD1}), we obtain
\begin{equation}\label{CaseIIb}
\begin{aligned}
\cosh l_{12}
&=1+\frac{2|z_1-z_2|^2}{(1-z_1^2)(1-|z_2|^2)}\\
&=\frac{(1+z_1^2)(1+|z_2|^2)-2z_1(z_2+\bar{z}_2)}{(1-z_1^2)(1-|z_2|^2)}\\
&=\frac{(1+xy)(1+|z_2|^2)-(x+y)(z_2+\bar{z}_2)}{\sqrt{1-x^2}\sqrt{1-y^2}(1-|z_2|^2)},
\end{aligned}
\end{equation}
where the last line follows from (\ref{Eq: F42}) and (\ref{Eq: F43}).
Combining (\ref{CaseIIa}) and (\ref{CaseIIb}) yields
\begin{equation}\label{Eq: F51}
2|z_2|^2+2xy-(x+y)(z_2+\bar{z}_2)
=\eta^h_{12}e^{u_2}(y-x)(1-|z_2|^2).
\end{equation}
By combining (\ref{Eq: F50}), (\ref{Eq: F51}), and the relationship $u_2=U_2+\ln\frac{2}{1-|z_2|^2}$, we conclude $\eta_{12}^h=\eta_{12}^e$.

\textbf{Case III}:
Assume that $\varepsilon_1=\varepsilon_2=1$.
See Figure \ref{figure4} (c).
The notations for the circle $C_1$ are consistent with those in Case 2.
Let $C_2$ denote the circle with hyperbolic center $v_2$.
The Euclidean and hyperbolic radii of $C_2$ are denoted by $R_2$ and $r_2$, respectively.
Let $z_2$ denote a complex number representing the position of $v_2$.
The intersection points of $C_2$ with the ray $Oz_2$ are denoted by $a$ and $b$.
Observe that the complex numbers $a$, $b$, and $z_2$ have the same argument, and $0<|a|<|z_2|<|b|<1$.

From (\ref{Eq: length-E3}), it follows that
\begin{align*}
\left|\frac{a+b}{2}-\frac{x+y}{2}\right|^2
&=L_{12}^2
=R_1^2+R_2^2+2R_1R_2\eta^e_{12}\\
&=\left(\frac{y-x}{2}\right)^2+\left|\frac{b-a}{2}\right|^2
+2\left(\frac{y-x}{2}\right)\left|\frac{b-a}{2}\right|\eta^e_{12}.
\end{align*}
This implies that
\begin{equation}\label{Eq: F15}
xy+|a||b|-\frac{1}{4}(x+y)(a+\bar{a}+b+\bar{b})
=\frac{1}{2}(y-x)|b-a|\eta_{12}^e.
\end{equation}
Note that in Case 2, the calculations involve real numbers $x$, $y$, and $z_1$.
We can replace the real numbers with complex numbers without difficulty.
By (\ref{Eq: HD2}), we obtain
\begin{equation}\label{Eq: etacheck i1j1}
\sinh r_2=\frac{|b-a|}{\sqrt{1-|a|^2}\sqrt{1-|b|^2}}
\quad \text{and} \quad
\cosh r_2=\frac{1-|a||b|}{\sqrt{1-|a|^2}\sqrt{1-|b|^2}}.
\end{equation}
From (\ref{Eq: length-H}), it follows that
\begin{equation}\label{Eq: F9}
\begin{aligned}
\cosh l_{12}
&=\cosh r_1\cosh r_2+\eta_{12}^h\sinh r_1\sinh r_2\\
&=\frac{(1-xy)(1-|a||b|)+\eta^h_{12}(y-x)|b-a|}{\sqrt{1-x^2}\sqrt{1-y^2}\sqrt{1-|a|^2}\sqrt{1-|b|^2}},
\end{aligned}
\end{equation}
where the second line follows from (\ref{Eq: etacheck i1j1}) and (\ref{Eq: F45}).
Furthermore,
\begin{equation*}
\frac{|z_2-a|}{\sqrt{1-|a|^2}\sqrt{1-|z_2|^2}}
=\sinh\frac{d(a,z_2)}{2}
=\sinh\frac{d(b,z_2)}{2}
=\frac{|b-z_2|}{\sqrt{1-|z_2|^2}\sqrt{1-|b|^2}}.
\end{equation*}
This implies that
\begin{equation}\label{Eq: F13}
z_2=\frac{a+b}{1+|a||b|+\sqrt{1-|a|^2}\sqrt{1-|b|^2}}.
\end{equation}
Consequently,
\begin{equation}\label{Eq: F14}
\begin{aligned}
1+|z_2|^2&=\frac{2(1+|a||b|)}{1+|a||b|+\sqrt{1-|a|^2}\sqrt{1-|b|^2}},\\
1-|z_2|^2&=\frac{2\sqrt{1-|a|^2}\sqrt{1-|b|^2}}{1+|a||b|+\sqrt{1-|a|^2}\sqrt{1-|b|^2}}.
\end{aligned}
\end{equation}
From (\ref{Eq: HD1}), it follows that
\begin{equation}\label{Eq: F12}
\begin{aligned}
\cosh l_{12}
&=1+\frac{2|z_1-z_2|^2}{(1-z_1^2)(1-|z_2|^2)}\\
&=\frac{(1+z_1^2)(1+|z_2|^2)-2z_1(z_2+\bar{z}_2)}{(1-z_1^2)(1-|z_2|^2)}\\
&=\frac{(1+xy)(1+|a||b|)-\frac{1}{2}(x+y)(a+\bar{a}+b+\bar{b})}{\sqrt{1-x^2}\sqrt{1-y^2}\sqrt{1-|a|^2}\sqrt{1-|b|^2}},
\end{aligned}
\end{equation}
where the last line follows from (\ref{Eq: F42}), (\ref{Eq: F43}), (\ref{Eq: F13}), and (\ref{Eq: F14}).
By comparing (\ref{Eq: F9}) and (\ref{Eq: F12}), we obtain
\begin{equation}\label{Eq: F16}
2(xy+|a||b|)-\frac{1}{2}(x+y)(a+\bar{a}+b+\bar{b})
=\eta^h_{12}(y-x)|b-a|.
\end{equation}
By combining (\ref{Eq: F15}) and (\ref{Eq: F16}), we conclude that $\eta_{12}^h=\eta_{12}^e$.

Note that in the preceding calculations, the Euclidean center of the circle $C_j$ associated with the vertex $v_j$ is replaced by its hyperbolic center.
In conjunction with the fact that the center of the circle at the vertex $v_0$ coincides with the origin,
it follows that the angle $\angle v_1v_0v_2$ is invariant computed in the PL metric $L$ or the PH metric $l$.
See Figure \ref{figure4}.

Based on the preceding analysis, we have the following lemma.

\begin{lemma}\label{Lem: PH-PL}
Let $l$ be the PH metric with hyperbolic label $u$ defined on $(P_n, \mathcal{T}, \varepsilon, \eta)$ with the interior vertex $v_0$ located at the origin $O$.
Let $L$ be the induced PL metric with Euclidean label $U$.
Then the following statements hold:
\begin{description}
\item[(a)] if $u$ and $U$ satisfy Convention \ref{Convention}, then $\eta^h=\eta^e$;
\item[(b)] the interior angle $\angle v_iv_0v_{i+1}$ in the PL metric $L$ is the same as that in the PH metric $l$ for all $i=1,\dots,n$;
\item[(c)] if $\eta^h=\eta^e$, then $l$ is weighted Delaunay if and only if $L$ is weighted Delaunay.
\end{description}
\end{lemma}



\subsection{Proof of the hyperbolic discrete maximal principle and the discrete Schwarz-Ahlfors lemma}\label{Sec: key2}

Recall that the discrete conformal factor in the Euclidean case is defined as $w^e=\tilde{U}-U$,
while in the hyperbolic case, it is defined as $w^h=\tilde{u}-u$.

Consider a geodesic segment $e=v_iv_j$ in $\mathbb{D}$.
Under a Euclidean or hyperbolic discrete conformal transformation, a new geodesic segment $\tilde{e}=\tilde{v}_i\tilde{v}_j$ is obtained in $\mathbb{D}$.
Without loss of generality, suppose that both $v_i$ and $\tilde{v}_i$ coincide with the origin.

If $\varepsilon_i=0$, then by Convention \ref{Convention} (a), it follows that
\begin{equation}\label{wie0}
w^h_i=\tilde{u}_i-u_i
=\tilde{U}_i+\ln 2-U_i-\ln 2
=w^e_i.
\end{equation}
If $\varepsilon_i=1$, then by Convention \ref{Convention} (b), it follows that
\begin{equation}\label{wie1}
w^h_i=\tilde{u}_i-u_i
=\tilde{U}_i-U_i
=w^e_i.
\end{equation}
If $\varepsilon_j=0$, then by Convention \ref{Convention} (c), it follows that
\begin{equation}\label{Eq: F52}
w^h_j=\tilde{u}_j-u_j
=\tilde{U}_j+\ln \frac{2}{1-|\tilde{z}_j|^2}
-U_j-\ln \frac{2}{1-|z_j|^2}=w^e_j+\ln\frac{1-|z_j|^2}{1-|\tilde{z}_j|^2},
\end{equation}
where $z_j$ and $\tilde{z}_j$ are complex numbers representing the positions of $v_j$ and $\tilde{v}_j$, respectively.
Note that the relationship (\ref{Eq: F52}) was first derived by Dai-Wu \cite[Lemma 2.3]{Dai-Wu}.
If $\varepsilon_j=1$, then by Convention \ref{Convention} (d), it follows that
\begin{equation}\label{Eq: F54}
\begin{aligned}
w^h_j
&=\tilde{u}_j-u_j\\
&=\tilde{U}_j+\ln \frac{2}{1-\tilde{x}\tilde{y}+\sqrt{(1-\tilde{x}^2)(1-\tilde{y}^2)}}
-U_j-\ln \frac{2}{1-xy+\sqrt{(1-x^2)(1-y^2)}}\\
&=w^e_j+\ln\frac{1-xy+\sqrt{(1-x^2)(1-y^2)}}{1-\tilde{x}\tilde{y}+\sqrt{(1-\tilde{x}^2)(1-\tilde{y}^2)}},
\end{aligned}
\end{equation}
where $x$ and $y$ (respectively, $\tilde{x}$ and $\tilde{y}$) denote the minimal and maximal Euclidean distances from the points of $C_j$ with hyperbolic center $v_j$ (respectively, $\tilde{C}_j$ with hyperbolic center $\tilde{v}_j$), to the origin $O$, respectively.
For a visual representation of the positions of the circles $C_j$ and $\tilde{C}_j$, please refer to Figure \ref{figure4} (c).

\begin{lemma}\label{Lem: relation}
Let $l$ be the PH metric with hyperbolic label $u$ defined on $(P_n, \mathcal{T}, \varepsilon, \eta)$ with the interior vertex $v_0$ located at the origin $O$.
Let $L$ be the induced PL metric with Euclidean label $U$.
Under the similarity transformation $z \to e^{\mu} z$ in the plane, $L$ induces a PL metric $\tilde{L}$ with Euclidean label $\tilde{U}$, which in turn induces a PH metric $\tilde{l}$ with hyperbolic label $\tilde{u}$.
If both $u,U$ and $\tilde{u},\tilde{U}$ satisfy Convention \ref{Convention},
then $l$ and $\tilde{l}$ are hyperbolic discrete conformal equivalent.
Moreover, 
\begin{equation}
w^h_0=w_0^e\equiv\mu, 
\end{equation}
and
\begin{equation}\label{Eq: F28}
w^h_j=\mu+\ln \frac{1-|z_j|^2}{1-e^{2\mu}|z_j|^2}
\end{equation}
for $\varepsilon_j=0$, and
\begin{equation}\label{Eq: F55}
w^h_j=\mu+\ln \frac{1-xy+\sqrt{(1-x^2)(1-y^2)}}{1-e^{2\mu}xy +\sqrt{(1-e^{2\mu}x^2)(1-e^{2\mu}y^2)}}
\end{equation}
for $\varepsilon_j=1$.
Here $v_j$ is adjacent to $v_0$, and $x$ and $y$ denote the minimal and maximal Euclidean distances, respectively, from the points of the circle $C_j$ (with hyperbolic center $v_j$) to the origin $O$.
\end{lemma}
\proof
By Lemma \ref{Lem: PH-PL} (a), we have $\eta^e=\eta^h$ and $\tilde{\eta}^e=\tilde{\eta}^h$.
Note that under a similarity transformation in the plane, the PL metrics $L$ and $\tilde{L}$ are Euclidean discrete conformal equivalent, i.e., $\eta^e=\tilde{\eta}^e$. 
Hence, $\eta^h=\tilde{\eta}^h$.
This implies that $l$ and $\tilde{l}$ are hyperbolic discrete conformal equivalent.
Furthermore, the similarity transformation of the PL metrics $L$ and $\tilde{L}$ implies that $w^e\equiv\mu$.
Consequently, $w_0^h=w_0^e\equiv\mu$ follows from (\ref{wie0}) and (\ref{wie1}).
The formulas (\ref{Eq: F28}) and (\ref{Eq: F55}) are derived from (\ref{Eq: F52}) and (\ref{Eq: F54}), respectively.
\qed

\begin{corollary}\label{Cor: key}
Under the same assumptions as in Lemma \ref{Lem: relation}, if $\mu<0$, then
\begin{equation*}
w_0^h>w_j^h,\ \forall j\sim 0.
\end{equation*}
\end{corollary}
\proof
If $\mu<0$, then
\begin{equation*}
\ln \frac{1-|z_j|^2}{1-e^{2\mu}|z_j|^2}<0 \quad \text{and} \quad
\ln \frac{1-xy+\sqrt{(1-x^2)(1-y^2)}}{1-e^{2\mu}xy +\sqrt{(1-e^{2\mu}x^2)(1-e^{2\mu}y^2)}}<0.
\end{equation*}
Therefore, by (\ref{Eq: F28}) and (\ref{Eq: F55}), it follows that $w_j^h<\mu=w_0^h$.
\qed

We now prove the hyperbolic maximum principle as follows.
\begin{theorem}\label{Thm: MP-H2}
Let $(\varepsilon,\eta)$ be a  pair of regular weights on $(P_n,\mathcal{T})$.
Suppose two PH metrics $l$ and $\bar{l}$ are hyperbolic discrete conformal equivalent, i.e., $\bar{l}=w^h\ast l$, and satisfy
\begin{enumerate}
\item[(a)] both $l$ and $\bar{l}$ are weighted Delaunay,
\item[(b)] the combinatorial curvatures $K_0(l)$ and $K_0(\bar{l})$ at the vertex $v_0$ satisfy $K_0(l) \geq K_0(\bar{l})$,
\item[(c)] all circles and points associated with hyperbolic labels $f$ and $\bar{f}$ of $l$ and $\bar{l}$ are contained in $\mathbb{D}$.
\end{enumerate}
Then the maximum of $w^h$, i.e. $\max_{j\in \{0, 1, \cdots, n\}}w^h_j=\max_{j\in \{0, 1, \cdots, n\}}(\bar{u}_j-u_j)$, if $>0$, is never achieved at $v_0$.
\end{theorem}

The condition (c) ensures the hyperbolic radii of circles are finite. 

\proof
We prove this theorem by contradiction.
Assume that
\begin{equation}\label{assumption}
w^h_0=\bar{u}_0-u_0
=\max_{j\in \{0, 1, \cdots, n\}}[\bar{u}_j-u_j]>0.
\end{equation}
By M\"{o}bius transformations, we can assume that $v_0=O$.
Then $l$ and $\bar{l}$ with hyperbolic labels $u$ and $\bar{u}$ induce two PL metrics $L$ and $\bar{L}$, respectively.
We choose the Euclidean labels $U$ and $\bar{U}$ of $L$ and $\bar{L}$ such that both $u,U$ and $\bar{u},\bar{U}$ satisfy Convention \ref{Convention}.
From Lemma \ref{Lem: PH-PL} (a), it follows that $\eta^h=\eta^e$ and $\bar{\eta}^h=\bar{\eta}^e$.
Since $\eta^h=\bar{\eta}^h$ by assumption, we have $\eta^e=\bar{\eta}^e$.
Therefore, $L$ and $\bar{L}$ are Euclidean discrete conformal equivalent.
Moreover, $L$ and $\bar{L}$ are weighted Delaunay by Lemma \ref{Lem: PH-PL} (c).
Set $\lambda=\frac{e^{u_0}}{e^{\bar{u}_0}}$. 
From (\ref{assumption}), it follows that
\begin{equation}\label{Eq: F18}
\lambda=\frac{e^{u_0}}{e^{\bar{u}_0}}<1.
\end{equation}
Applying a similarity transformation $z \to \lambda z$ to $\bar{L}$ in the plane,
a PL metric $\bar{L}^\lambda$ with Euclidean label $\bar{U}^\lambda$ and the corresponding PH metric $\bar{l}^\lambda$ are obtained.
Note that $\bar{U}^\lambda$ is determined by $\bar{U}$.
Since the similarity transformation preserves the weighted Delaunay condition,
$\bar{L}^\lambda$ is weighted Delaunay.
Choose the hyperbolic label $\bar{u}^\lambda$ of $\bar{l}^\lambda$ such that $\bar{u}^\lambda,\bar{U}^\lambda$ satisfy Convention \ref{Convention}.
By Lemma \ref{Lem: PH-PL} (a) (c), $\bar{l}^{\lambda}$ is weighted Delaunay.
By Lemma \ref{Lem: relation}, $\bar{l}^{\lambda}$ is hyperbolic discrete conformal equivalent to $\bar{l}$.
Hence, $l$ and $\bar{l}^{\lambda}$ are hyperbolic discrete conformal equivalent, and $L$ and $\bar L^\lambda$ are Euclidean discrete conformal equivalent.

Applying Corollary \ref{Cor: key} to $\bar{l}$ and $\bar{l}^\lambda$, for $j\sim 0$, yields
\begin{equation*}
\bar u^\lambda_j-\bar u_j<\bar u_0^\lambda-\bar u_0 \quad \Leftrightarrow \quad
\bar u_0-\bar u_j<\bar u_0^\lambda-\bar u_j^\lambda.
\end{equation*}
Combining this result with the assumption that $\bar u_0-u_0\geq\bar u_j-u_j$ in (\ref{assumption}), we obtain
\begin{equation}\label{Eq: F27}
u_0-u_j\leq \bar u_0-\bar u_j<\bar u_0^\lambda-\bar u_j^\lambda.
\end{equation}
Since $\bar u_0^\lambda-\bar u_0=\ln \lambda$ by Lemma \ref{Lem: relation},
it follows that $e^{\bar u_0^\lambda}=\lambda e^{\bar u_0}=e^{u_0}$ by (\ref{Eq: F18}).
Consequently, $\bar u_0^\lambda=u_0$ and $\bar{f}_0^\lambda=f_0$.
From (\ref{Eq: F27}), it follows that $u_j>\bar u_j^\lambda$ and $f_j>\bar f_j^\lambda$.
Let $l_{0j}$ and $\bar l^\lambda_{0j}$ be the hyperbolic lengths of the edges from $v_j$ to the origin.
From (\ref{Eq: SC1}) and (\ref{Eq: length-H}), it follows that $l_{0j}>\bar l^\lambda_{0j}$ (here we use $\eta_{0j}=\bar{\eta}^{\lambda}_{0j}$). Let $L_{0j}$ and $\bar L^\lambda_{0j}$ be the Euclidean lengths of the edges from $v_j$ to the origin. Then
\begin{equation*}
L^2_{0j}
=\varepsilon_0e^{2U_0}+\varepsilon_je^{2U_j}
+2\eta^e_{0j}e^{U_0+U_j}
\quad \text{and} \quad
(\bar L^\lambda_{0j})^2
=\varepsilon_0e^{2\bar{U}_0^\lambda}+ \varepsilon_je^{2\bar{U}_j^\lambda}+
2\eta^e_{0j}e^{\bar{U}_0^\lambda+\bar{U}_j^\lambda}.
\end{equation*}
Since $\bar u_0^\lambda-u_0=\bar U_0^\lambda-U_0$ by Lemma \ref{Lem: relation},
it follows that $\bar u_0^\lambda=u_0$ implies $\bar U_0^\lambda=U_0$.
If $\varepsilon_j=0$, using the relation $L=\tanh \frac{l}{2}$ between the Euclidean and hyperbolic lengths,
then $l_{0j}>\bar l^\lambda_{0j}$ implies $L_{0j}>\bar L^\lambda_{0j}$.
Consequently, $U_j>\bar{U}_j^\lambda$ by (\ref{Eq: SC1}).
If $\varepsilon_j=1$, then from (\ref{Eq: length-H}), it follows that
\begin{equation*}
\cosh l_{0j}=\sqrt{1+\varepsilon_0e^{2f_0}}
\cosh r_j+\eta^h_{0j}e^{f_0}\sinh r_j.
\end{equation*}
Therefore,
\begin{align*}
e^{U_j}
&=\frac{1}{2}(\tanh \frac{l_{0j}+r_j}{2}-\tanh \frac{l_{0j}-r_j}{2})\\
&=\frac{\sinh r_j}{\cosh l_{0j}+\cosh r_j}\\
&=\frac{1}{(\sqrt{1+\varepsilon_0e^{2f_0}}+1)
\coth r_j+\eta^h_{0j}e^{f_0}}.
\end{align*}
Similarly,
\begin{equation*}
e^{\bar U_j^\lambda}
=\frac{1}{(\sqrt{1+\varepsilon_0e^{2\bar{f}^\lambda_0}}+1)
\coth \bar{r}^\lambda_j+\eta^h_{0j}e^{\bar{f}^\lambda_0}}.
\end{equation*}
Note that $u_j>\bar{u}^\lambda_j$ implies $r_j>\bar{r}^\lambda_j$.
Combining this with $f_0=\bar{f}_0^\lambda$ yields $U_j>\bar U_j^\lambda$.

In a 1-ring neighborhood of $v_0=O$, it holds that $U_0=\bar U_0^\lambda$ and $U_j>\bar{U}_j^\lambda$ for all $j\sim 0$.
This implies that the maximum of $e^{\bar{U}_j^\lambda}/e^{U_j}$ is attained at the interior vertex $v_0$.
Note that by Lemma \ref{Lem: PH-PL} (b), the hyperbolic angle at $O$ of a hyperbolic triangle coincides with the Euclidean angle at $O$ of the corresponding Euclidean triangle.
Consequently, the hyperbolic combinatorial curvature at $v_0$ equals the Euclidean combinatorial curvature at $v_0$.
By the assumption, it follows that $K_0(L)\geq K_0(\bar{L})$.
Since Euclidean angles remain invariant under a similarity transformation,
it follows that $K_0(\bar{L})=K_0(\bar{L}^\lambda)$.
Therefore, $K_0(L)\geq K_0(\bar{L}^\lambda)$.
By Theorem \ref{Thm: MP-E2}, $e^{\bar{U}_j^\lambda}/e^{U_j}=e^{\bar{U}_0^\lambda}/e^{U_0}$ for all $j\sim 0$.
Hence, $U_j=\bar U_j^\lambda$.
This leads to a contradiction.
\qed
\\

\noindent\textbf{Proof of Theorem \ref{Thm: MP-H}:}
The part (1) follows directly from Theorem \ref{Thm: MP-H2}.
We derive the part (2) by substituting $w^h_0$ with $-w^h_0$ in the part (1).
This operation is reasonable by the assumption that $K_0(l)\leq K_0(\bar{l})$.
\qed

\begin{remark}
Theorem \ref{Thm: MP-H} generalizes the hyperbolic maximal principles in Lemma 2.2 of He \cite{He} and Theorem 2.7 of Dai-Wu \cite{Dai-Wu}. Notably, our result does not require $P_n$ to be embeddable into the hyperbolic plane $\mathbb{D}$, i.e., $K_0(l)= K_0(\bar{l})= 0$.
\end{remark}

As an application of Theorem \ref{Thm: MP-H2}, we have the following discrete Schwarz-Ahlfors lemma.

\begin{theorem}[Discrete Schwarz-Ahlfors lemma]\label{Thm: DSL}
Let $(\varepsilon,\eta)$ be a regular pair of weights on $(M,\mathcal{T},\varepsilon,\eta)$, where $M\subseteq\mathbb{D}$ is a compact set with non-empty boundary.
Suppose $l$ and $\bar{l}$ are two weighted Delaunay PH metrics with hyperbolic labels $f$ and $\bar{f}$, respectively, satisfying all circles and points with respect to $f,\bar f$ are contained in $\mathbb{D}$.
\begin{description}
\item[(a)] If the combinatorial curvatures $K(l)\geq K(\bar{l})$ for all interior vertices, and $w^h \leq 0$ holds for every boundary vertex, then $w^h \leq 0$ holds for all vertices.
\item[(b)] If $K(l)\leq K(\bar{l})$ for all interior vertices, and $w^h \geq 0$ holds for every boundary vertex, then $w^h \geq 0$ holds for all vertices.
\end{description}
\end{theorem}
\proof
The first part follows from Theorem \ref{Thm: MP-H} (a). By contradiction, suppose there exists an interior vertex $v_i$ with $w^h_i>0$. Then $w^h$ achieve its maximum at some interior vertex.  Without loss of generality, we assume $w^h_i=\max_{j} w^h_j >0$.
By applying Theorem \ref{Thm: MP-H} (a) to the 1-ring neighborhood of $v_i$,
we deduce that there exists a vertex $j$ adjacent to $i$ such that $w^h_j>w^h_i$.
This contradicts the maximality of $w^h_i$.
The second part follows analogously from Theorem \ref{Thm: MP-H} (b).
\qed

Theorem \ref{Thm: DSL2} is a direct corollary of Theorem \ref{Thm: DSL}. Also, we have the following result related to the rigidity of PH metrics in discrete conformal structures.  
\begin{corollary}\label{rigidity of HDCF}
Under the same conditions as in Theorem \ref{Thm: DSL}, if $K(l)\equiv K(\bar{l})$ for all interior vertices, and $w^h\equiv 0$ holds for every boundary vertex, then $w^h\equiv0$ holds for all vertices.
\end{corollary}

In \cite{Xu 1}, there is another way to prove Corollary \ref{rigidity of HDCF} by constructing convex energy functions. Furthermore, using the method in \cite{Xu 1} to prove the rigidity, one does not
need to assume that the two PH metrics $l$ and $\bar{l}$ are weighted Delaunay. Please refer to \cite{Xu 1} for more details.

In the following, we explain why we refer to Theorem \ref{Thm: DSL} as \textit{Discrete Schwarz-Ahlfors Lemma}. The classical Schwarz-Ahlfors Lemma in \cite{Ah} could be stated as follows.
\begin{lemma}[Schwarz-Ahlfors Lemma \cite{Ah}]
\label{Lem: SAL}
Let $\mathbb{D}$ be the unit disk in the complex plane, and let $\rho(z)|dz|$ be the Poincar\'{e} metric in $\mathbb{D}$. Suppose $ds=\sigma(z)|dz|$ is a $C^2$ Riemannian metric on $\mathbb{D}$ with Gaussian curvature $K\leq -1$.
Then $\rho(z)\geq \sigma(z)$ for any $z\in \mathbb{D}$.
\end{lemma}

Using Ahlfors's original proof of Schwarz-Ahlfors lemma, one can prove the following generalized Schwarz-Ahlfors lemma.
\begin{lemma}[Generalized Schwarz-Ahlfors Lemma]
\label{Lem: GSAL}
Let $D$ be a domain in the complex plane, and let $ds=\sigma_j|\mathrm{d}z|$ be the metrics on $D$ for $j=1,2$, where $\sigma_j>0$ and $\sigma_j \in C^2$. 
Assume that the Gaussian curvatures satisfy $K(\sigma_1)\geq K(\sigma_2)$ and $K(\sigma_2)<0$,
and that for every boundary point $\xi \in \partial D$, $\liminf_{z \to \xi} \frac{\sigma_1(z)}{\sigma_2(z)} \geq 1$.
Then $\sigma_1(z) \geq \sigma_2(z)$ for all $z\in D$.
\end{lemma}

\proof
This should be a well-known result, but we did not find a proof of the result in the literature.
A related work on complete Riemannian surfaces is Troyanov \cite{Tr}.
For completeness, we give a proof here. The idea of the proof comes from the original work of Ahlfors \cite{Ah}.

Set $F(z)=\ln\sigma_1(z)-\ln\sigma_2(z)$. We just need to prove that $F(z)\geq 0$ in $D$.
By contradiction, suppose that there is a point $z_0\in D$ such that $F(z_0)<0$. Then $\sigma_1(z_0) < \sigma_2(z_0)$.

By the assumption, we have $\liminf_{z\to\xi}F(z)\geq 0$ for all $\xi\in\partial D$.
Hence, for any each $\xi\in\partial D$,
there exists an open neighborhood $U_\xi$ of $\xi\in\partial D$ such that $F(z)> \frac{1}{2}F(z_0)$ for all $z\in U_\xi\cap D$.
By finite covering theorem, there exists a finite number of open sets $U_1, \ldots, U_n$ such that $K=D \setminus \left(\bigcup_{\alpha=1}^n U_\alpha \right)$
is a non-empty compact set.
As a result, $F(z)$ attains its minimum in $D$. Without loss of generality, we can assume
$$F(z_0)=\inf_{z\in D}F(z)<0.$$
Then $\Delta F(z_0)\geq0$.

On the other hand, at $z_0$, we have
\begin{align*}
\Delta F|_{z_0}
=&(\Delta\ln\sigma_1-\Delta\ln\sigma_2)|_{z_0}\\
=&[-K(\sigma_1)\sigma^2_1+K(\sigma_2)\sigma^2_2]|_{z_0}\\
\leq& -K(\sigma_2)(\sigma^2_1-\sigma^2_2)|_{z_0}\\
<&0,
\end{align*}
where the formula $K(\sigma)=-\sigma^{-2}\Delta \ln \sigma$ for the Gaussian curvature is used in the second line and the condition $K(\sigma_2)<0$ is used in the third line.
This contradicts $\Delta F(z_0)\geq0$.
\qed

It is straightforward to see that Theorem \ref{Thm: DSL2} is a discretization of the generalized Schwarz-Ahlfors Lemma, i.e. Theorem \ref{Lem: GSAL}.
Notably, Theorem \ref{Thm: DSL2} generalizes prior discrete Schwarz lemmas in \cite{He, BS2, BS3,Van1994}. 
Please refer to \cite{Rodin1, Rodin2, HS1993, HS1996} and others for different versions of discrete Schwarz lemma.

\section{Infinite rigidity of small Delaunay triangulations of the plane}
\label{Sec: IR}

In this section, we prove Theorem \ref{Thm: IR-H} as an application of Theorem \ref{Thm: MP-H} in the special case of  $\varepsilon\equiv0$.
Consequently, the Euclidean and hyperbolic discrete conformal structures reduce to Euclidean and hyperbolic vertex scalings, respectively.
Recall that if two PL metrics $L$ and $\tilde L$ are Euclidean discrete conformal equivalent,
then
\begin{equation}\label{Eq: F20}
\tilde L_{ij}=L_{ij}e^{\frac{1}{2}(w^e_i+w^e_j)}.
\end{equation}
If two PH metrics $l$ and $\tilde l$ are hyperbolic discrete conformal equivalent, then
\begin{equation}\label{Eq: DCE-H2}
\sinh\frac{\tilde l_{ij}}{2}
=e^{\frac{1}{2}(w^h_i+w^h_j)}\sinh \frac{l_{ij}}{2}.
\end{equation}

The following lemma is a special case of Lemma \ref{Lem: relation} when $\epsilon\equiv 0$. Recall $\mathcal{T}$ is an infinite geodesic triangulation of $\mathbb{D}$ in Theorem \ref{Thm: IR-H}.

\begin{lemma}\label{Lem: DCE}
Let $\mathcal{T}_0=(V_0,E_0,F_0)$ be a subcomplex of $\mathcal{T}$. Let $z_i$ be the coordinate of $v_i\in V_0$ and $\tilde z_i$ be the image of $z_i$ under similarity transformation $z\rightarrow e^\mu z$ in the plane. Let $L$ and $\tilde L$ be the PL metrics on $\mathcal{T}_0$ determined by $\{z_i\}$ and $\{\tilde z_i\}$. Let $l$ and $\tilde l$ be the PH metrics on $\mathcal{T}_0$ determined by $\{z_i\}$ and $\{\tilde z_i\}$. 
Then $l$ and $\tilde{l}$ are discrete conformal equivalent, i.e., $\tilde{l}=w^h*l$, with
\begin{equation}\label{Eq: F53}
w^h_i=\mu+\ln \frac{1-|z_i|^2}{1-|\tilde{z}_i|^2}.
\end{equation}
\end{lemma}

\begin{lemma}\label{Lem: dist}
The following statements hold.
\begin{description}
\item[(a)]
For all $x\in[0, 1]$, the inequality $\arcsin x\leq 2x$ holds, and 
\begin{equation*}
\frac{x}{2}\leq \sin x \leq x\leq \sinh x\leq e^x - 1\leq 2x.
\end{equation*}
\item[(b)]
If the hyperbolic distance $d(z_1, z_2)<1$, then
\begin{equation*}
\frac{|z_2-z_1|}{1-|z_2|} \leq 2d(z_1, z_2).
\end{equation*}

\item[(c)]
If $|z_1|\leq |z_2|$, then
\begin{equation*}
\frac{1}{2}d(z_1, z_2) \leq \frac{|z_1-z_2|}{1-|z_2|}.
\end{equation*}
\end{description}
\end{lemma}
\begin{proof}
\textbf{(a):} The proof of this part is straightforward and thus omitted here.

\textbf{(b):}
Let $C\subseteq\mathbb{H}^2$ be the hyperbolic circle containing $z_1$ and $z_2$ such that the diameter of $C$ is the geodesic segment from $z_1$ to $z_2$. 
Apply a rotational isometry centered at the hyperbolic center of $C$ to map $z_1$ and $z_2$ to $z_1'$ and $z_2'$, so that $z_1'$ and $z_2'$ lie on the $x$-axis with $z_1', z_2'\in (-1,1)$, $z_2'>0$, and $|z_1'|\leq z_2'$.
It is straightforward to check that 
$$\frac{|z_2-z_1|}{1-|z_2|} \leq \frac{|z_2'-z_1'|}{1-|z_2'|},$$
and $d(z_1',z_2') = d(z_1, z_2)$. Hence, it suffices to show the case of $z_1 = z_1'$ and $z_2 = z_2'$. 

The cross-ratio formula for the hyperbolic distance yields
\begin{equation*}
d(z_1,z_2)
=\ln \frac{(1-z_1)(1+z_2)}{(1-z_2)(1+z_1)}
\geq \ln\frac{1-z_1}{1-z_2}
=\ln \left(1 +\frac{|z_2-z_1|}{1-|z_2|}\right).
\end{equation*}
Furthermore, the inequality $e^x-1\leq 2x$ implies that
\begin{equation*}
\frac{|z_2-z_1|}{1-|z_2|}\leq e^{d(z_1,z_2)}-1
\leq 2d(z_1, z_2).
\end{equation*}

\textbf{(c):}
Let $\gamma(t)=tz_2+(1-t)z_1$ be the parametrized curve that connects $z_1$ and $z_2$. 
Clearly, $|\gamma(t)|\leq t|z_2|+(1-t)|z_1|\leq |z_2|$ holds for all $t\in[0,1]$.
Therefore,
\begin{equation*}
d(z_2,z_1)
\leq l(\gamma)
\leq \int_0^1\frac{2|\gamma'(t)|}{1-|\gamma(t)|^2}dt \leq \frac{2|z_2-z_1|}{1-|z_2|}.
\end{equation*}
\end{proof}

For any three points not on the same geodesic in $\mathbb{D}$, there exist both a Euclidean triangle and a hyperbolic triangle.
When the triangle is sufficiently small under the hyperbolic metric, replacing the hyperbolic geodesic with a Euclidean geodesic results in only minimal changes to the triangle's angles as described below. 
\begin{lemma}\label{angle}
Let $\gamma$ and $\gamma'$ be the hyperbolic and Euclidean geodesic segments in $\mathbb{D}$ connecting $z_1$ and $z_2$, respectively.
If $d(z_1, z_2)\leq 1$, then the angle of intersection between $\gamma$ and $\gamma'$ is less than $2d(z_1, z_2)$.
\end{lemma}
\begin{proof}
If $\gamma$ extends to a diameter of $\mathbb{D}$, then $\gamma'$ coincides with $\gamma$, and the angle of intersection is zero.
Otherwise, let $c$ be the center of the circle containing $\gamma$ in $\mathbb{D}$, where the circle is orthogonal to the unit circle.
Since $|c|>1$, it follows that $|c-z_1|>1-|z_1|$.
Furthermore, combining Lemma \ref{Lem: dist} (b) with the assumption that $d(z_1,z_2)\leq 1$ yields
\begin{equation*}
\frac{|z_2-z_1|}{2|c-z_1|}
\leq \frac{|z_2-z_1|}{2(1-|z_1|)}
\leq d(z_1, z_2)\leq 1.
\end{equation*}
Using the inequality $\arcsin x \leq 2x$ in Lemma \ref{Lem: dist}(a) when $x\in [0, 1]$, the angle of intersection between $\gamma$ and $\gamma'$ is
\begin{equation*}
\frac{1}{2}\angle z_1cz_2
=\arcsin \frac{|z_2-z_1|}{2|c-z_1|}
\leq \frac{|z_2-z_1|}{|c-z_1|}
<\frac{|z_2-z_1|}{1-|z_1|}
\leq 2d(z_1,z_2).
\end{equation*}
\end{proof}

\noindent\textbf{Proof of Theorem \ref{Thm: IR-H}:}
By the assumption, there exists a discrete conformal factor $w^h$ such that $\bar{l}=w^h*l$.
It suffices to show that $w_i^h\geq 0$ for all $v_i\in V$.
By symmetry, we have $w^h\equiv0$, which implies that $\bar{l}$ and $l$ are the same PH metrics.

Otherwise, assume that $w^h_0<0$ for some vertex $v_0$.
Without loss of generality, we assume that $v_0=O$ for both triangulations $l$ and $\bar l$.
We perform the following transformation to $\bar{l}$ while preserving $l$.

\textbf{Step 1:}
Let $z_i$ be the complex number representing the position of $v_i$ in $\bar{l}$.
Then $\{z_i\}$ forms a set of points in $\mathbb{D}$.
Replace the hyperbolic geodesic segment connecting $z_i$ and $z_j$ with the Euclidean segment connecting $z_i$ and $z_j$.
In general, the PH metric $\bar{l}$ does not necessarily induce a PL metric $\bar{L}$.
However, under the assumptions in Theorem \ref{Thm: IR-H}, the PH metric $\bar{l}$ induces a PL metric $\bar{L}$ in this way.
By the assumption that all inner angles of all triangles in $\bar{l}$ are at least $\delta$, we have $\delta<\frac{\pi}{3}$.
Since the maximum length of the edge in $\bar l$ is assumed to be less than $\delta^3/8192$, by Lemma \ref{angle}, the angles of the Euclidean triangles in $\bar{L}$ are at least $\delta/2$, and the orientation of the triangles is preserved. 
Thus, $\bar{L}$ is a PL metric induced by $\bar{l}$.

\textbf{Step 2:}
Apply a similarity transformation $z\to e^\mu z$ to $\bar L$, where $0<\mu<-w^h_0$ is chosen sufficiently small to obtain a new PL metric $\tilde L$.
Under this transformation, the Euclidean angles remain unchanged.
Let $\tilde z_i$ be the complex number in $\tilde{L}$ representing the position of $v_i$.
Note that some $\tilde z_i$ may lie outside the unit disk due to the scaling.
Recall that $d(z_1,z_2)=\infty$ if either $z_1$ or $z_2$ lies outside $\mathbb{D}$.
However, we are not concerned with points outside $\mathbb{D}$.
Instead, we define a subcomplex $T_0$ of $\mathcal{T}$ with the following properties:
\begin{enumerate}
\item[(1)] all edges $v_iv_j$ in $E_0=E(T_0)$ satisfy $d(\tilde z_i, \tilde z_j)\leq \delta/16$,
\item[(2)] all vertices in $V_0=V(T_0)$ are the set of vertices intersecting some edges in $E_0$,
\item[(3)] all faces in $F_0=F(T_0)$ are the set of triangles whose edges are all in $E_0$.
\end{enumerate}
Note that $T_0$ is a finite complex, since the set $\{z_i\}$ is discrete in $\mathbb{D}$ with its limit set in $\partial \mathbb{D}$. Since $\mu$ is chosen to be sufficiently small,  it is straightforward to check that $v_0\in V_0$. 

\textbf{Step 3:}
Replace the Euclidean geodesic segments in $\tilde{L}$ with hyperbolic geodesic segments to obtain a hyperbolic metric $\tilde l$.
By Lemma \ref{angle} again, $\tilde l$ is a PH metric on $T_0$ with angles bounded by $\delta/4$.
This completes the construction.

By Lemma \ref{Lem: DCE}, the PH metrics $l$ and $\tilde l$ are discrete conformal equivalent, i.e., $\tilde l =\tilde w^h*l$ with
\begin{equation*}
\tilde w^h_i
=w^h_i+\mu+\ln\frac{1-|z_i|^2}{1-|\tilde z_i|^2}.
\end{equation*}
This follows from the fact that $\bar l=w^h*l$ by definition and $\tilde{l}=w* \bar l$ with $w_i=\mu_i+\ln\frac{1-|z_i|^2}{1-|\tilde z_i|^2}$ by (\ref{Eq: F53}).

By $v_0=O$ and the definition of $\mu$, it follows that $\tilde w^h_0=w^h_0+\mu<0$.
Let $\tilde w^h_i=\min_{j\in V_0} \tilde w^h_j$.
By the hyperbolic maximum principle in Theorem \ref{Thm: MP-H}, $\tilde w^h_i$ is attained at a boundary vertex $v_i\in V_0$.
Furthermore, there exists $v_k$ adjacent to $v_i$ such that $v_k\notin V_0$.
Otherwise, all the vertices adjacent to $v_i$ are contained in $V_0$.
Applying the hyperbolic maximum principle in Theorem \ref{Thm: MP-H}
to the 1-ring neighborhood of $v_i$, we have $\tilde w^h_i>\min_{j, v_iv_j\in E(\mathcal{T})} \tilde w^h_j$, which contradicts $\tilde w^h_i=\min_{j\in V_0} \tilde w^h_j$.
As a result, there exists a triangle $\triangle v_iv_jv_k$ such that $v_iv_j\in E_0$ and $v_iv_k\not\in E_0$.
Therefore, by the definition of $E_0$, we have $\tilde l_{ij}\leq \delta/16$ and $\tilde l_{ik}>\delta/16$.

Since $\tilde l_{ij}\leq \delta/16$, we can bound $\tilde l_{ik}$ and $\tilde l_{jk}$ using the Euclidean sine law as follows.
By Lemma \ref{Lem: dist} (b), the inequality $\tilde l_{ij}\leq \delta/16$ implies that
\begin{equation}\label{Eq: F21}
\frac{|\tilde z_i-\tilde z_j|}{1-|\tilde z_j|} 
\leq 2\tilde l_{ij} \leq\frac{\delta}{8}.
\end{equation}
Note that the inner angles of the Euclidean triangle $\triangle_E v_iv_jv_k$ are at least $\delta/2$.
By the Euclidean sine law, we have  
\begin{equation}
\label{Eq71}
\frac{|\tilde z_k-\tilde z_j|}{1-|\tilde z_j|}
\leq \frac{|\tilde z_i-\tilde z_j|}{1-|\tilde z_j|}\frac{1}{\sin(\delta/2)}
\leq \frac{\delta}{8\sin(\delta/2)}
\leq \frac{\delta}{8(\delta/2)/2}
=\frac{1}{2},
\end{equation}
where the third inequality follows from $\sin x\geq x/2$ in Lemma \ref{Lem: dist} (a).
This, combined with the inequality $|\tilde z_k-\tilde z_j|+|\tilde z_j|\geq |\tilde z_k|$, implies that
\begin{equation*}
1-|\tilde z_k|
\geq 1-|\tilde z_j|-|\tilde z_j-\tilde z_k|
\geq 1-|\tilde z_j|-\frac{1}{2}(1-|\tilde z_j|)
=\frac{1}{2}(1-|\tilde z_j|).
\end{equation*}
From Lemma \ref{Lem: dist} (c), it follows that
\begin{equation*}
\frac{1}{2}\tilde l_{jk}
\leq \max \left\{\frac{|\tilde z_k - \tilde z_j|}{1- |\tilde z_k|}, \frac{|\tilde z_k - \tilde z_j|}{1-|\tilde z_j|}\right\}
\leq \frac{|\tilde z_k - \tilde z_j|}{\frac{1}{2}(1-|\tilde z_j|)}.
\end{equation*}
Therefore, 
\begin{equation}
\label{Eq72}
\frac{1}{2}\tilde l_{jk}
\leq\frac{2|\tilde z_k-\tilde z_j|}{1-|\tilde z_j|}
\leq\frac{2|\tilde z_j-\tilde z_i|}{(1-|\tilde z_j|)\sin(\delta/2)}
\leq\frac{8}{\delta}\frac{|\tilde z_i-\tilde z_j|}{1-|\tilde z_j|}
\leq\frac{8}{\delta}\cdot 2\tilde l_{ij}
\leq 1,
\end{equation}
where the second inequality follows from the Euclidean sine law, the third inequality follows from $\sin x\geq x/2$ in Lemma \ref{Lem: dist} (a), and the forth inequality follows from (\ref{Eq: F21}).
Consequently, $\tilde l_{jk}\leq (32/\delta)\tilde l_{ij}\leq 2$. 
Similarly, by exchanging the index $j$ with $i$ in (\ref{Eq: F21}) and replacing the index $j$ by $i$ in (\ref{Eq71}), we can show that $\tilde l_{ik}\leq 2$ by repeating the computation above.

In summary, the edges in the hyperbolic triangle $\triangle_H v_iv_jv_k$ satisfy the following bounds: 
$\tilde l_{ij}\leq \delta/16$, $\tilde l_{jk}\leq 2$, and $\delta/16< \tilde l_{ik}\leq 2$.
Then, from the definition of $\tilde w^h_i$ in (\ref{Eq: DCE-H2}), it follows that
\begin{equation*}
e^{\tilde w^h_i}
=\frac{(\sinh \frac{\tilde l_{ij}}{2}/\sinh \frac{l_{ij}}{2})(\sinh \frac{\tilde l_{ik}}{2}/\sinh \frac{l_{ik}}{2})}{(\sinh \frac{\tilde l_{jk}}{2}/\sinh \frac{l_{jk}}{2})}
=\frac{\sinh \frac{\tilde l_{ik}}{2}}{\sinh \frac{l_{ij}}{2}}\frac{\sinh\frac{\tilde l_{ij}}{2}}{\sinh \frac{\tilde l_{jk}}{2}}\frac{\sinh\frac{l_{jk}}{2}}{\sinh \frac{ l_{ik}}{2}}.
\end{equation*}
Therefore,
\begin{equation*}
e^{\tilde w^h_i}
\geq \frac{\tilde l_{ik}/2}{l_{ij}}\frac{\tilde l_{ij}/2}{\tilde l_{jk}}\frac{(\sinh l_{jk})/4}{\sinh l_{ik}/2}
=\frac{1}{8}\frac{\tilde l_{ik}}{l_{ij}}\frac{\tilde l_{ij}}{\tilde l_{jk}}\frac{\sinh l_{jk}}{\sinh l_{ik}}
\geq\frac{1}{8}\frac{\delta/16}{l_{ij}}\frac{\delta}{32}\sin\delta
\geq \frac{\delta^3}{8192l_{ij}}
\geq 1,
\end{equation*}
where the first inequality follows from $x\geq\sinh (x/2)\geq x/2\geq (\sinh x)/4$ in Lemma \ref{Lem: dist} (a) and $2\sinh (l_{ik}/2)\leq \sinh l_{ik}$,
the second inequality follows from $\tilde l_{ij}/\tilde l_{jk}\geq \delta/32$ by (\ref{Eq72}), $\tilde l_{ik}> \delta/16$, and $\sinh l_{jk}/\sinh l_{ik}\geq \sin\delta$ by the hyperbolic cosine law and the assumption,
the third inequality follows from $\sin \delta \geq \delta/2$,
and the last inequality follows from $l_{ij}\leq \delta^3/8192$.
This contradicts with that $\tilde w^h_i<\tilde w^h_0<0$.
\qed

\end{document}